\documentclass[11pt,a4paper]{article}

\usepackage{amsmath}
\usepackage{amssymb}

\numberwithin{equation}{section}

\usepackage{color}

\newcommand{\p}{\partial}
\newcommand{\nn}{\nonumber}

\newcommand{\C}{\mathbb{C}}
\newcommand{\A}{\mathcal{A}}
\newcommand{\hA}{\hat{\A}}
\newcommand{\F}{\mathcal{F}}
\newcommand{\hF}{\hat{\F}}
\newcommand{\V}{\mathcal{V}}
\renewcommand{\H}{\mathcal{H}}
\newcommand{\X}{\mathcal{X}}
\renewcommand{\S}{\mathcal{S}}
\newcommand{\Ker}{\mathrm{Ker}}
\newcommand{\ad}{\mathrm{ad}}
\newcommand{\tV}{\hat{\V}}
\newcommand{\tH}{\hat{\H}}
\newcommand{\tX}{\hat{\X}}
\newcommand{\tnabla}{\hat{\nabla}}

\newcommand{\al}{\alpha}

\newtheorem{dfn}{Definition}[section]
\newtheorem{lem}[dfn]{Lemma}
\newtheorem{prp}[dfn]{Proposition}
\newtheorem{thm}[dfn]{Theorem}
\newtheorem{rmk}[dfn]{Remark}
\newtheorem{cor}[dfn]{Corollary}
\newtheorem{emp}[dfn]{Example}

\newenvironment{prf}{\noindent \textit{Proof} \ }{\hfill $\Box$}
\newenvironment{prfn}[1]{\noindent \textit{Proof of #1} \ }{\hfill $\Box$}

\begin{document}

\title{{Bihamiltonian Cohomologies and Integrable Hierarchies II: the Tau Structures}}
\author{
Boris Dubrovin$^{*}$, Si-Qi Liu$^\dag$, Youjin Zhang$^\dag$\\
{\small ${}^*$ SISSA, via Bonomea 265, Trieste 34136, Italy}\\
{\small ${}^{\dag}$ Department of Mathematical Sciences, Tsinghua University,}\\
{\small Beijing 100084, P. R. China}}
\renewcommand{\thefootnote}{}
\footnotetext{Emails: dubrovin@sissa.it, liusq@tsinghua.edu.cn, youjin@tsinghua.edu.cn}
\renewcommand{\thefootnote}{\arabic{footnote}} 
\date{}\maketitle

\begin{abstract}
Starting from a so-called flat exact semisimple bihamiltonian structures of hydrodynamic type, we arrive at a Frobenius manifold structure and a tau structure for the associated principal hierarchy. We then classify the deformations of the principal hierarchy which possess tau structures.
\end{abstract}

{\small
\noindent\textbf{Mathematics Subject Classification (2010).} Primary 37K10; Secondary 53D45.}


\section{Introduction}
The class of bihamiltonian integrable hierarchies which possess hydrodynamic limits plays an important role in the study of Gromov--Witten invariants,
2D topological field theory, and other research fields of mathematical physics. In \cite{DZ-NF} the first- and third-named authors of the present paper
initiated a program of classifying deformations of bihamiltonian integrable hierarchies of hydrodynamic type under the so-called Miura type transformations.
They introduced the notion of bihamiltonian cohomologies of a bihamiltonian structure and converted the classification problem into the computation of
these cohomology groups. The first two bihamiltonian cohomologies for semisimple bihamiltonian structures of hydrodynamic type were calculated in
\cite{DLZ-1, LZ-1}, and it was proved that the infinitesimal deformations of a semisimple bihamiltonian structure of hydrodynamic type are parametrized
by a set of smooth functions of one variable. For a given deformation of a semisimple bihamiltonian structure of hydrodynamic type these functions
$c_1(u^1),\dots, c_n(u^n)$ can be calculated by an explicit formula represented in terms of the canonical coordinates $u^1,\dots, u^n$ of the bihamiltonian
structure. These functions are invariant under the Miura type transformations, due to this reason they are called the central invariants of the deformed
bihamiltonian structure.

In \cite{BCIH-I}, the second- and third-named author of the present paper continued the study of the above mentioned classification problem.
They reformulated the notion of infinite dimensional Hamiltonian structures in terms of the infinite jet space of a super manifold, and provided
a framework of infinite dimensional Hamiltonian structures which is convenient for the study of properties of Hamiltonian and bihamiltonian
cohomologies.  One of the main results which is crucial for the computation of bihamiltonian cohomologies is given by Lemma 3.7 of \cite{BCIH-I}.
It reduces the computation of the bihamiltonian cohomologies to the computations of cohomology groups of a bicomplex on the space of differential
polynomials, instead of on the space of local functionals. Based on this result, they computed the third bihamiltonian cohomology group of the
bihamiltonian structure of the dispersionless KdV hierarchy, and showed that any infinitesimal deformation of this bihamiltonian structure can be
extended to a full deformation.

In \cite{CPS-2},  Carlet, Posthuma and Shadrin completed the computation of the third bihamiltonian cohomology group for a general semisimple
bihamiltonian structure of hydrodynamic type based on the results of \cite{BCIH-I}.  Their result confirms the validity of the conjecture of \cite{BCIH-I}
that any infinitesimal deformation of a semisimple bihamiltonian structures of hydrodynamic type can be extended to a full deformation, i.e. for any given
smooth functions $c_i(u^i)\ (i=1, \dots, n)$, there exists a deformation of the corresponding semisimple bihamiltonian structure of hydrodynamic type such
that its central invariants are given by $c_i(u^i)\ (i=1, \dots, n)$. 

This paper is a continuation of \cite{BCIH-I}. We are to give a detailed study of properties of the integrable hierarchies associated with a special class of
semisimple bihamiltonian structures of hydrodynamic type and their deformations, which are called \emph{flat exact semisimple bihamiltonian structures
of hydrodynamic type}. One of their most important properties is the existence of tau structures for the associated integrable hierarchies and their
deformations with constant central invariants.

For a hierarchy of Hamiltionian evolutionary PDEs, a tau structure is a suitable choice of the densities of the Hamiltonians satisfying certain conditions
which enables one to define a function, called the tau function, for solutions of the hierarchy of evolutionary PDEs, as it is defined in \cite{DZ-NF}.
The notion of tau functions was first introduced by M.~Sato \cite{Sato} for solutions to the KP equation and  by Jimbo, Miwa and Ueno for a class of
monodromy preserving deformation equations of linear ODEs with rational coefficients \cite{JMU-1, JM-2, JM-3} at the beginning of 80's of the last
century. It was also adopted to soliton equations that can be represented as equations of isospectral deformations of certain linear spectral problems
or as Hamiltonian systems, and has played crucial role in the study of relations of soliton equations with  infinite dimensional Lie algebras \cite{DKJM, KW},
and with the geometry of infinite dimensional Grassmannians \cite{SS, SW}. The importance of the notion of tau functions of soliton equations is manifested
by the discovery of the fact that the tau function of a particular solution of the KdV hierarchy is a partition function of 2D gravity, see \cite{Wi, Ko} for
details. In \cite{DZ-NF}, the first- and the third-named authors introduced the notion of tau structures for the class of bihamiltonian integrable hierarchies
possessing hydrodynamic limits, and constructed the so-called topological deformations of the principal hierarchy of a semisimple Frobenius manifold by
using properties of the associated tau functions. On the other hand, not all bihamiltonian integrable hierarchies possess tau structures.  In this paper we
introduce the notion of \emph{flat exact} bihamiltonian structure, and study the classification of the associated tau structures. It turns out that this notion
is an appropriate generalization of semisimple conformal Frobenius manifolds when considering the associated integrable hierarchies and their tau
structures. One can further consider the deformations of a flat exact semisimple bihamiltonian structure of hydrodynamic type which possess tau structures.
It is known that the central invariants of such deformations must be constant \cite{yz}. We show that deformations with constant central invariants of
a flat exact semisimple bihamiltonian structure of hydrodynamic type indeed possess tau structures, and we also give a classification theorem for the
associated tau structures.

The paper is arranged as follows. In Sec.\,\ref{sec-1} we introduce the notion of flat exact  semisimple bihamiltonian structures of hydrodynamic type and
present the main results. In Sec.\,\ref{sec-2} we study the relations between flat exact  semisimple bihamiltonian structures of hydrodynamic type and
semisimple Frobenius manifolds, and give a proof of Theorem \ref{mainthm00}. In Sec.\,\ref{sec-3} we construct the principal hierarchy for a flat exact
semisimple bihamiltonian structures of hydrodynamic type and show the existence of a tau structure. In Sec.\,\ref{sec-4} we consider properties of
deformations of the principal hierarchies which possess tau structures and the Galilean symmetry, and then in Sec.\,\ref{sec-5} we prove the existence of
deformations of the principal hierarchy of a flat exact bihamiltonian structures of hydrodynamic type, which are bihamiltonian integrable hierarchies
possessing tau structures and the Galilean symmetry, and we prove Theorem \ref{main-thm}.  Sec.\,\ref{sec-7} is a conclusion. In the Appendix, we prove some properties of
semi-Hamiltonian  integrable hierarchies, some of which are used in the proof of the uniqueness theorem given in Sec.\,\ref{sec-4}.

\section{Some notions and the main results }\label{sec-1}

The class of systems of hydrodynamic type on the infinite jet space of an $n$-dimensional manifold $M$  consists of systems of $n$ first order quasilinear
partial differential equations (PDEs)  
\begin{equation}\label{sht00}
v^\alpha_t = \sum_{\beta=1}^n A_\beta^\alpha(v) v^\beta_x, \quad \alpha=1, \dots, n, \quad v=\left( v^1, \dots, v^n\right)\in M.
\end{equation}
Here $A^\alpha_\beta(v)$ is a section of the bundle $TM \otimes T^*M$. For the subclass of Hamiltonian systems of hydrodynamic type the r.h.s. of
\eqref{sht00} admits a representation
\begin{equation}\label{hsht00}
v^\alpha_t = P^{\alpha\beta} \frac{\partial h(v)}{\partial v^\beta}.
\end{equation}
Here the smooth function $h(v)$ is the density of the Hamiltonian
\[H=\int h(v)\, dx\]
and
\begin{equation}\label{pbht00}
P^{\alpha\beta}=g^{\alpha\beta}(v) \partial_x +\Gamma^{\alpha\beta}_\gamma(v) v^\gamma_x
\end{equation}
is the operator of a \emph{Poisson bracket of hydrodynamic type}. As it was observed in \cite{DN83} such operators satisfying the \emph{nondegeneracy condition}
\begin{equation}
\det \left( g^{\alpha\beta}(v)\right)\neq 0 \label{nondegeneracy}
\end{equation}
correspond to flat metrics (Riemannian or pseudo-Riemannian)
\[ds^2 =g_{\alpha\beta}(v) dv^\alpha dv^\beta\]
on the manifold $M$. Namely,
\[g^{\alpha\beta}(v)=\left( g_{\alpha\beta}(v)\right)^{-1}\]
is the corresponding inner product on $T^*M$, the coefficients $\Gamma^{\alpha\beta}_\gamma(v)$ are the contravariant components of the Levi-Civita
connection for the metric. In the present paper it will be assumed that all Poisson brackets of hydrodynamic type satisfy the nondegeneracy condition
\eqref{nondegeneracy}. 

A bihamiltonian structure of hydrodynamic type is a pair $(P_1, P_2)$ of operators of the form \eqref{pbht00} such that an arbitrary linear combination
$\lambda_1 P_1+\lambda_2 P_2$ is again the operator of a Poisson bracket. They correspond to pairs of flat metrics $g_1^{\alpha\beta}(v)$,
$g^{\alpha\beta}_2(v)$ on $M$ satisfying certain compatibility condition (see below for the details). The bihamiltonian structure of hydrodynamic type is
called \emph{semisimple} if the roots $u^1(v)$, \dots, $u^n(v)$ of the characteristic equation
\begin{equation}\label{char00}
\det \left( g^{\alpha\beta}_2(v)-u \,g^{\alpha\beta}_1(v)\right)=0
\end{equation}
are pairwise distinct and are not contant for a generic point $v\in M$. 
According to Ferapontov's theorem \cite{Fera}, these roots can serve as local coordinates
of the manifold $M$, which are called the canonical coordinates of the bihamiltonian structure $(P_1, P_2)$.
We assume in this paper that $D$ is a sufficiently small domain on $M$ such that $(u^1, \dots, u^n)$ is the local 
coordinate system on $D$. In the canonical coordinates the two metrics have diagonal forms
\begin{equation}\label{zh-10-30}
g_1^{ij}(u)=f^i(u)\delta^{ij}, \quad g_2^{ij}(u)=u^i\,f^i(u)\delta^{ij}.
\end{equation}
We will need to use the notion of rotation coefficients of the metric $g_1$ 
which are defined by the following formulae:
\begin{equation}\label{zh-11-1}
\gamma_{ij}(u)=\frac{1}{2\sqrt{f_i f_j}} \frac{\p f_i}{\p u^j},\quad i\ne j
\end{equation}
with $f_i=\frac{1}{f^i}$. We also define $\gamma_{ii}=0$.
\begin{dfn}[cf. \cite{DZ-NF}]\label{zh-11-2}
The semisimple bihamiltonian structure $(P_1, P_2)$ is called reducible at $u\in M$
if there exists a partition of the set $\{1, 2, \dots, n\}$ into the union of 
two nonempty nonintersecting sets $I$ and $J$ such that
\[
 \gamma_{ij}(u)= 0,\quad \forall i\in I,\ \forall j\in J.\]
$(P_1, P_2)$ is called irreducible on a certain domain $D\subset M$, if it is not reducible at any point $u\in D$.
\end{dfn}

The main goal of the present paper is to introduce tau-functions of bihamiltonian systems of hydrodynamic type and of their dispersive deformations.
This will be done under the following additional assumption.

\begin{dfn}\label{zh-10-31}
The bihamiltonian structure $(P_1, P_2)$ of hydrodynamic type is called \emph{exact} if there exists a vector field $Z\in Vect\left( M\right)$ such that
\begin{equation}\label{cond-exact}
[Z, P_1]=0, \quad [Z, P_2]=P_1.
\end{equation}
Here $[\ \,, \ ]$ is the infinite-dimensional analogue of the Schouten--Nijenhuis bracket
(see the next section and \cite{BCIH-I} for details of the definition).
It is called \emph{flat exact} if the vector field $Z$ is flat with respect to the metric associated with the Hamiltonian structure $P_1$. 
\end{dfn}

\begin{emp} \label{exam-frob} Let $\left( M, \cdot\,, \eta, e, E\right)$ be a Frobenius manifold. Then the pair of metrics 
\begin{equation}\label{pencil-frob}
\begin{aligned}
& g_1^{\alpha\beta}(v)=\langle dv^\alpha, dv^\beta\rangle= \eta^{\alpha\beta},\\
& g_2^{\alpha\beta}(v)=(dv^\alpha, dv^\beta)=i_E \left( dv^\alpha\cdot dv^\beta\right)=:g^{\alpha\beta}(v)
\end{aligned}
\end{equation}
on $T^*M$ defines a flat exact bihamiltonian structure with $Z=e$, see \cite{Du-1} for the details. For a semisimple Frobenius manifold the resulting
bihamiltonian structure will be semisimple. Roots of the characteristic equation \eqref{char00} coincide with the canonical coordinates on the Frobenius
manifold.
\end{emp}

More bihamiltonian structures can be obtained from those of Example \ref{exam-frob} by a Legendre-type transformation
\cite{Du-1, XZ}
\begin{equation}\label{legen01}
\hat v_\alpha =b^\gamma\frac{\p^2 F(v)}{\p v^\gamma\p v^\alpha}, \quad \hat v^\alpha=\eta^{\alpha\beta}\hat v_\beta.
\end{equation}
Here $F(v)$ is the potential of the Frobenius manifold and $b=b^\gamma\frac{\p}{\p v^\gamma}$ is a flat invertible vector field on it.
The new metrics on $T^*M$ by definition have the \emph{same} Gram matrices in the new coordinates
\begin{equation}\label{legen02}
\langle d\hat v^\alpha, d\hat v^\beta\rangle =\eta^{\alpha\beta}, \quad \left(  d\hat v^\alpha, d\hat v^\beta\right) =g^{\alpha\beta}(v).
\end{equation}
Recall that applying the transformation \eqref{legen01} to $F(v)$ one obtains a new solution $\hat F(\hat v)$ to the WDVV associativity equations defined from
\begin{equation}\label{legen03}
\frac{\p^2 \hat{F}(\hat{v})}{\p\hat{v}^\al\p\hat{v}^\beta}=\frac{\p^2 F(v)}{\p{v}^\al\p{v}^\beta}.
\end{equation}
The new unit vector field is given by
\begin{equation}\label{legen04}
\hat e =b^\gamma\frac{\p}{\p\hat v^\gamma}.
\end{equation}
The new solution to the WDVV associativity equations defines on $M$ another Frobenius manifold structure if the vector $b=b^\gamma\frac{\p}{\p v^\gamma}$ satisfies
\[\left[ b, E\right] = \lambda\cdot b\]
for some $\lambda\in\mathbb C$. Otherwise the quasihomogeneity axiom does not hold true.

\begin{thm} \label{mainthm00}
For an arbitrary Frobenius manifold $M$ the pair of flat metrics obtained from \eqref{pencil-frob} by a transformation of the form
\eqref{legen01}--\eqref{legen02} defines on $M$ a flat exact bihamiltonian structure of hydrodynamic type. Conversely, any irreducible
flat exact semisimple bihamiltonian structure of hydrodynamic type can be obtained in this way.
\end{thm}

Now we can describe a tau-symmetric bihamiltonian hierarchy associated with a flat exact semisimple bihamiltonian structure $(P_1, P_2; Z)$ of
hydrodynamic type. Let us choose a system of flat coordinates $\left( v^1, \dots, v^n\right)$ for the first metric. So the operator $P_1$ has the form
\[P_1^{\alpha\beta}=\eta^{\alpha\beta}\frac{\p}{\p x}\]
for a constant symmetric nondegenerate matrix $\eta^{\alpha\beta}=g_1^{\alpha\beta}$. 
It is convenient to normalize the choice of flat coordinates by the requirement
\[Z=\frac{\p}{\p v^1}.\]
We are looking for an infinite family of systems of first order quasilinear evolutionary PDEs of the form \eqref{sht00}
satisfying certain additional conditions. The systems of the form \eqref{sht00} will be labeled by pairs of indices $(\alpha, p)$, $\alpha=1, \dots, n$, $p\geq 0$.
Same labels will be used for the corresponding time variables $t=t^{\alpha, p}$. The conditions to be imposed are as follows.

1. All the systems under consideration are \emph{bihamiltonian} PDEs w.r.t. $(P_1, P_2)$. This implies pairwise commutativity of the flows \cite{DLZ-1}
\begin{equation}\label{comm00}
\frac{\p}{\p t^{\alpha,p}}\frac{\p v^\gamma}{\p t^{\beta,q}}=\frac{\p}{\p t^{\beta,q}}\frac{\p v^\gamma}{\p t^{\alpha,p}}.
\end{equation}

2. Denote
\begin{equation}\label{hamilt00}
H_{\alpha,p}=\int h_{\alpha, p}(v)\, dx
\end{equation}
the Hamiltonian of the $(\alpha,p)$-flow with respect to the first Poisson bracket,
\begin{equation}\label{hamilt01}
\frac{\p v^\gamma}{\p t^{\alpha,p}}=\eta^{\gamma\lambda}\frac{\p}{\p x} \frac{\delta H_{\alpha,p}}{\delta v^\lambda(x)} \equiv \eta^{\gamma\lambda}\frac{\p}{\p x} \frac{\p h_{\alpha,p}(v)}{\p v^\lambda}.
\end{equation}
The Hamiltonian densities satisfy the following \emph{recursion}\footnote{This recursion acts in the opposite direction with respect to the
bihamiltonian one - see eq. \eqref{bi-recursion00} below.}
\begin{equation}\label{recur00}
\frac{\p}{\p v^1} h_{\alpha, p}(v) = h_{\alpha,p-1}(v), \quad \alpha=1, \dots, n, \quad p\geq 0
\end{equation}
(recall that $\frac{\p}{\p v^1} =Z$) where we denote
\begin{equation}\label{casi00}
h_{\alpha, -1}(v)=v_\alpha\equiv\eta_{\alpha\beta}v^\beta, \quad \alpha=1, \dots, n.
\end{equation}
Observe that the functionals $H_{\alpha,-1}=\int h_{\alpha,-1}(v)\, dx$ span the space of Casimirs of the first Poisson bracket.

3. Normalization
\begin{equation}\label{normalize01}
\frac{\p}{\p t^{1,0}}=\frac{\p}{\p x}.
\end{equation}

\begin{prp}\label{prp-25}
Integrable hierarchies of the above form satisfy the \emph{tau-symmetry} condition
\begin{equation}\label{tau-sym00}
\frac{\p h_{\alpha, p-1}}{\p t^{\beta,q}}= \frac{\p h_{\beta,q-1}}{\p t^{\alpha, p}},\quad \forall ~ \alpha, \, \beta=1, \dots, n, \quad \forall~ p, \, q\geq 0.
\end{equation}
Moreover, this integrable hierarchy is invariant with respect to the  Galilean symmetry
\begin{align}
&\frac{\p v}{\p s} =Z(v)+\sum_{p\geq 1} t^{\alpha, p}\frac{\p v}{\p t^{\alpha, p-1}},\label{galileo}\\
&\left[ \frac{\p}{\p s}, \frac{\p}{\p t^{\alpha,p}}\right]=0, \quad \forall \alpha=1, \dots, n, \quad p\geq 0.\nn
\end{align}
\end{prp}

\begin{dfn}\label{dfn-cali}
A choice of the Hamiltonian densities $h_{\alpha, p}(v)$, $\alpha=1, \dots, n$, $p\geq -1$ satisfying the above conditions is called a \emph{calibration}
of the flat exact bihamiltonian structure $(P_1, P_2; Z)$ of hydrodynamic type. The integrable hierarchy \eqref{hamilt01} is called the
\emph{principal hierarchy} of $(P_1, P_2; Z)$ associated with the given calibration.
\end{dfn}

\begin{emp} Let $\left( M, \,\cdot\,, \eta, e, E\right)$ be a Frobenius manifold. Denote
\[\left(\theta_1(v; z), \dots, \theta_n(v; z)\right) z^\mu z^R\]
with
\begin{align}
&\theta_\alpha(v; z) =\sum_{p=0}^\infty \theta_{\alpha,p}(v) z^p, \quad \alpha=1, \dots, n\label{levelt00}
\end{align}
a Levelt basis of deformed flat coordinates \cite{Du-3}. Here the matrices 
\[\mu={\rm diag}(\mu_1, \dots, \mu_n), \quad R=R_1+\dots,\quad  [\mu, R_k] =k\, R_k\] 
constitute a part of the spectrum of the Frobenius manifold, see details in \cite{Du-3}. Then 
\begin{equation}\label{levelt01}
h_{\alpha,p}(v) =\frac{\p \theta_{\alpha,p+2}(v)}{\p v^1}, \quad \alpha=1, \dots, n, \quad p\geq -1
\end{equation}
is a calibration of the flat exact bihamiltonian structure associated with the metrics \eqref{pencil-frob} on the Frobenius manifold. In this case the family of
pairwise commuting bihamiltonian PDEs \eqref{hamilt01} is called the \emph{principle hierarchy} associated with the Frobenius manifold. With this choice
of the calibration the Hamiltonians \eqref{hamilt00}, \eqref{levelt01} satisfy the bihamilonian recursion relation
\begin{equation}\label{bi-recursion00}
\{ \, .\,, H_{\beta,q-1}\}_2=
(q+\frac12+\mu_\beta)\{ \, . \,, H_{\beta,q}\}_1+
\sum_{k=1}^{q-1} (R_{q-k})^\al_\gamma 
\{ \, . \,, H_{\beta,k}\}_1.
\end{equation}
Other calibrations can be obtained by taking constant linear combinations and shifts
\begin{align}
& \tilde\theta_\alpha(v; z) =\theta_\beta(v; z) C_\alpha^\beta(z)+\theta_{\alpha}^0(z), \quad \alpha =1, \dots, n \label{basis-change00} \\
& C(z)=\left( C_\alpha^\beta(z)\right) =\mathbf{1}+C_1 z+C_2 z^2+\dots, \quad C^T(-z) C(z)=\mathbf{1} \nn\\
& \theta_\alpha^0(z)=\sum_{p\geq 0}\theta_{\alpha,p}^0 z^p, \quad \theta_{\alpha,p}^0\in \mathbb C. \nn
\end{align}
\end{emp}

For the flat exact bihamiltonian structure obtained from \eqref{pencil-frob} by a Legendre-type transformation \eqref{legen01}--\eqref{legen04}
one can choose a calibration by introducing functions $\hat \theta_{\alpha,p}\left( \hat v\right)$ defined by
\begin{equation}\label{legen05}
\frac{\p \hat\theta_\alpha\left(\hat v; z\right)}{\p \hat v^\beta}=\frac{\p \theta_\alpha(v; z)}{\p v^\beta},\quad \forall\, \alpha, \, \beta=1, \dots, n.
\end{equation}
Remarkably in this case the new Hamiltonians satisfy the \emph{same} bihamiltonian recursion \eqref{bi-recursion00}.
Other calibrations can be obtained by transformations of the form \eqref{basis-change00}.

\begin{prp}
For a flat exact bihamiltonian structure of hydrodynamic type obtained from a Frobenius manifold by a Legendre-type transformation \eqref{legen01}--\eqref{legen04} the construction \eqref{legen05} and \eqref{levelt01} defines a calibration. Any calibration can be obtained in this way up to the transformation \eqref{basis-change00} .
\end{prp}

The properties of a calibration, in particular the tau-symmetry property 
\eqref{tau-sym00}, of a flat exact semisimple bihamiltonian structure of hydrodynamic type $(P_1, P_2; Z)$ enable us to define a tau structure and tau functions for it and the associated principal hierarchy \eqref{hamilt01}, see Definitions \ref{zh-12-2} and \ref{zh-01-22f} in Section \ref{sec-3}. 
One of the main purposes of the present paper is to study the existence and properties of tau structures for deformations of the bihamiltonian structure $(P_1, P_2; Z)$ and the principal hierarchy.
Let $c_i(u^i)\ (i=1, \dots, n)$ be a collection of arbitrary smooth functions,  Carlet, Posthuma, and Shadrin showed that there exists a deformation
$(\tilde{P}_1, \tilde{P}_2)$ of $(P_1, P_2)$ such that its central invariants are given by $c_i(u^i)\ (i=1, \dots, n)$ \cite{CPS-3}.
By using the triviality of the second bihamiltonian cohomology,  one can show that there also exists a unique deformation of the principal hierarchy
of $(P_1, P_2)$ such that all its members are bihamiltonian vector fields of $(\tilde{P}_1, \tilde{P}_2)$ (see Sec.\,\ref{sec-5}). The deformed
integrable hierarchy usually does not possess a tau structure unless the central invariants are constant (first observed in \cite{yz}). On the other hand,
it is shown by Falqui and Lorenzoni in \cite{FL} that, if $c_i(u^i)\ (i=1, \dots, n)$ are constants, one can choose the representative $(\tilde{P}_1, \tilde{P}_2)$
such that they still satisfy the exactness condition, that is
\[[Z, \tilde{P}_1]=0, \quad [Z, \tilde{P}_2]=\tilde{P}_1.\]
With such a pair $(\tilde{P}_1, \tilde{P}_2)$ in hand, we can ask the following questions:
\begin{enumerate}
\item Does the deformed integrable hierarchy have tau structures?
\item If it does, how many of them? 
\end{enumerate}
The following theorem is the main result of the present paper, which answers the above questions.

\begin{thm}\label{main-thm}
Let $(P_1, P_2; Z)$ be a flat exact semisimple bihamiltonian structure of hydrodynamic type which satisfies the irreducibility condition. We fix a calibration 
$\{h_{\al,p}\,|\,\al=1,\dots,n;\, p=0,1,2,\dots\}$ of the bihamiltonian structure 
$(P_1, P_2; Z)$.
Then the following statements hold true:
\begin{itemize}
\item[i)] For any deformation $(\tilde{P}_1, \tilde{P}_2; \tilde{Z})$ of $(P_1, P_2; Z)$ with constant central invariants, there exists a deformation 
$\{\tilde{h}_{\alpha, p}\}$ of the Hamiltonian densities $\{h_{\al, p}\}$ such that the corresponding Hamiltonian vector fields $\tilde{X}_{\alpha, p}$
yield a deformation of the principal hierarchy which is a bihamiltonian integrable hierarchy possessing a tau structure and the Galilean symmetry.
\item[ii)] Let $(\hat{P}_1, \hat{P}_2; \hat{Z})$ be another deformation of $(P_1, P_2; Z)$ 
with the same central invariants as $(\tilde{P}_1, \tilde{P}_2; \tilde{Z})$, and let $\{\hat{h}_{\alpha, p}\}$ be the corresponding tau-symmetric
deformation of the Hamiltonian densities, then the logarithm of the tau function for $\{\hat{h}_{\alpha, p}\}$ can be obtained from the one for
$\{\tilde{h}_{\alpha, p}\}$ by adding a differential polynomial.
\end{itemize}
\end{thm}

\section{Flat exact semisimple bihamitonian structures and Frobenius manifolds}\label{sec-2}

Let $M$ be a smooth manifold of dimension $n$. Denote by $\hat{M}$ the super manifold of
dimension $(n\mid n)$ obtained from the cotangent bundle of $M$ by reversing the parity of the 
fibers. Suppose $U$ is a local coordinate chart on $M$ with coordinates $(u^1, \dots, u^n)$,
then
\[\theta_i=\frac{\p}{\p u^i},\quad i=1, \dots, n\]
can be regarded as local coordinates on the corresponding local chart $\hat{U}$ on $\hat{M}$.
Note that $\theta_i$'s are super variables, they satisfy the skew-symmetric commutation law:
\[\theta_i\theta_j+\theta_j\theta_i=0.\]

Let $J^\infty(M)$ and $J^\infty(\hat{M})$ be the infinite jet space of $M$ and $\hat{M}$,
which is just the projective limits of the corresponding finite jet bundles. There is a natural local
chart $\hat{U}^\infty$ over $\hat{U}$ with local coordinates
\[\{u^{i,s}, \theta_i^s \mid i=1, \dots, n; s=0, 1, 2, \dots\}.\]
See \cite{BCIH-I} for more details. Denote by $\hA$ the spaces of differential polynomials on $\hat{M}$.
Locally, we can regard $\hA$ as
\[C^\infty(\hat{U})[[u^{i,s}, \theta_i^s \mid i=1, \dots, n; s=1, 2, \dots]].\]
The differential polynomial algebra $\A$ on $M$ can be defined similarly as a subalgebra of $\hA$.
There is a globally defined derivation on $J^\infty(\hat{M})$
\begin{equation}\label{zh-12-1}
\p=\sum_{i=1}^n\sum_{s\ge0}\left(u^{i,s+1}\frac{\p}{\p u^{i,s}}+\theta_i^{s+1}\frac{\p}{\p \theta_i^s}\right).
\end{equation}
Its cokernel $\hF=\hA/\p\hA$ is called the space of local functionals. Denote the projection $\hA\to\hF$ by $\int$. We can also define $\F=\A/\p\A$, whose elements
are called local functionals on $M$.

There are two useful degrees on $\hA$, which are called standard gradation
$$
\deg u^{i,s}=\deg \theta_i^s=s
$$ 
and super gradation
$$
\deg \theta_i^s=1, \quad \deg u^{i,s}=0 
$$
respectively:
\[\hA=\bigoplus_{d\ge 0}\hA_d=\bigoplus_{p\ge 0}\hA^p.\]
We denote $\hA^p_d=\hA_d\cap \hA^p$.
In particular, $\A=\hA^0$, $\A_d=\hA^0_d$. The derivation $\p$ has the property $\p (\hA^p_d) \subseteq \hA^p_{d+1}$, hence it induces the same degrees on $\hF$, so we also have the homogeneous components $\hF_d$, $\hF^p$, $\hF_d^p$,
and the ones for $\F=\hF^0$. The reader can refer to \cite{BCIH-I} for details of the definitions of these notations.

There is a graded Lie algebra structure on $\hF$, whose bracket operation is given by
\[[P, Q]=\int\left(\frac{\delta P}{\delta \theta_i}\frac{\delta Q}{\delta u^i}+(-1)^p\frac{\delta P}{\delta u^i}\frac{\delta Q}{\delta \theta_i}\right),\]
where $P\in\hF^p$, $Q\in \hF^q$. This bracket is called the Schouten--Nijenhuis bracket on $J^\infty(M)$.

A Hamiltonian structure is defined as an elements $P\in\hF^2$ satisfying $[P, P]=0$. For example,
the operator \eqref{pbht00} corresponds to an element $P\in\hF^2_1$ of the form
\[P=\frac12\int\left(g^{ij}(u)\theta_i\theta_j^1+\Gamma^{ij}_{k}(u)u^{k,1}\theta_i\theta_j\right).\]
The fact that $P$ is a Hamiltonian operator is equivalent to the condition $[P, P]=0$.

A bihamiltonian structure of hydrodynamic type can be given by a pair of Hamiltonian structures of hydrodynamic type
$(P_1, P_2)$ satisfying the additional condition $[P_1, P_2]=0$. Denote by $g_1, g_2$ the flat metrics associated with the Hamiltonian structures $P_1, P_2$.
In what follows, we will assume that $(P_1, P_2)$ is semisimple with a fixed system of canonical coordinates $u^1, \dots, u^n$, in which the two flat metrics take the 
diagonal form \eqref{zh-10-30}, and the contravariant Christoffel coefficients of them have the following expressions respectively:
\begin{align}
\Gamma^{ij}_{k}&=\frac12\frac{\p f^i}{\p u^k}\delta^{ij}+
\frac12\frac{f^i}{f^j}\frac{\p f^j}{\p u^i}\delta^{jk}-\frac12\frac{f^j}{f^i}\frac{\p f^i}{\p u^j}\delta^{ik},\label{gamma-1}\\
\hat\Gamma^{ij}_{k}&=\frac12\frac{\p (u^i f^i)}{\p u^k}\delta^{ij}+
\frac12\frac{u^i f^i}{f^j}\frac{\p f^j}{\p u^i}\delta^{jk}-\frac12\frac{u^j f^j}{f^i}\frac{\p f^i}{\p u^j}\delta^{ik}\label{gamma-2}.
\end{align}
The diagonal entries $f^i$ satisfy certain non-linear differential equations which are equivalent to the flatness of $g_1$, $g_2$ and the condition
$[P_1, P_2]=0$. See the appendix of \cite{DLZ-1} for details. We denote by $\nabla$, $\tnabla$ the Levi-Civita connections of the metrics $g_1$, $g_2$
respectively.

We also assume henceforth that the semisimple bihamiltonian structure of hydrodynamic type $(P_1, P_2)$ 
is flat exact (see Definition \ref{zh-10-31}), and the corresponding vector field is given by $Z\in\hF^1$. We will denote this exact 
bihamiltonian structure by $(P_1, P_2; Z)$. 

\begin{lem} \label{lem-24}
If $Z\in\hF^1$ satisfies the condition \eqref{cond-exact}, then it has the following form:
\[Z=\int \left(\sum_{i=1}^n \theta_i\right)+X,\]
where $X$ is a bihamiltonian vector field of $(P_1, P_2)$.
\end{lem}
\begin{prf}
We first decompose $Z\in \hF^1$ into the sum of homogeneous components:
\[Z=Z_0+Z_1+Z_2+\cdots, \quad \mbox{where } Z_k\in\hF^1_k.\]
It is proved in \cite{FL} that $Z_0$ must take the form 
\begin{equation}
Z_0=\int \left(\sum_{i=1}^n \theta_i\right). \label{eq-Z0}
\end{equation}
Then $X:=Z-Z_0$ satisfies $[X, P_1]=[X, P_2]=0$, so it is a bihamiltonian vector field of $(P_1, P_2)$.
\end{prf}

The $X$-part of $Z$ does not affect anything, so it can be omitted safely.
Then $Z=Z_0$, and we call it the unit vector field of $(P_1, P_2)$.
According to the convention used in \cite{BCIH-I}, this $Z$ corresponds to a vector field on $M$ given by
\[D_Z=\sum_{i=1}^n \frac{\p}{\p u^i}\]
(see Definition 2.2 and Equation (2.5) of \cite{BCIH-I}).
It is also proved in \cite{FL} that if \eqref{cond-exact} holds true then 
\begin{equation}\label{w-10}
D_Z(f^i)=\sum_{k=1}^n \frac{\p f^i}{\p u^k}=0,\quad i=1, \dots, n. 
\end{equation}

Note that the flatness of the vector field $Z$ (or, equivalently, $D_Z$)  given in Definition \ref{zh-10-31} can be 
represented as
\begin{equation}\label{cond-flat}
 \nabla D_Z=0.
 \end{equation}

\begin{lem} \label{lem-25}
$D_Z$ is flat if and only if $f_i:=(f^i)^{-1}\ (i=1, \dots, n)$ satisfy the following Egoroff conditions:
\begin{equation}\label{jw-5-2}
\frac{\p f_i}{\p u^j}=\frac{\p f_j}{\p u^i}, \quad \forall\,1\le i,j \le n.
\end{equation}
\end{lem}
\begin{prf}
The components of $D_Z$ read $Z^j=1$, so we have
\begin{equation}\label{zh-07}
0=\nabla^i Z^j=g_1^{ik}\frac{\p Z^j}{\p u^k}-\Gamma^{ij}_{k} Z^k=-\sum_{k=1}^n \Gamma^{ij}_{k}.
\end{equation}
By using \eqref{gamma-1} and \eqref{w-10}, the lemma can be easily proved.
\end{prf}

The above lemma implies that, if $Z$ is flat, then $\gamma_{ij}=\gamma_{ji}$ (see \eqref{zh-11-1}).
In this case, the conditions that $(P_1, P_2)$ is a bihamiltonian structure
are equivalent to the following equations for $\gamma$ (see the appendix of \cite{DLZ-1}):
\begin{align}
&\frac{\p \gamma_{ij}}{\p u^k}=\gamma_{ik}\gamma_{jk}, \quad \mbox{for distinct } i, j, k, \label{gamma-cond-1}\\
&\sum_{k=1}^n \frac{\p \gamma_{ij}}{\p u^k}=0, \label{gamma-cond-2}\\
&\sum_{k=1}^n u^k \frac{\p \gamma_{ij}}{\p u^k}=-\gamma_{ij}. \label{gamma-cond-3}
\end{align}
The condition \eqref{gamma-cond-2} is actually $D_Z(\gamma_{ij})=0$. If we introduce the Euler vector field
\begin{equation}
E=\sum_{k=1}^n u^k \frac{\p}{\p u^k}, \label{euler-vf}
\end{equation}
then the condition \eqref{gamma-cond-3} is $E(\gamma_{ij})=-\gamma_{ij}$, that is, $\gamma_{ij}$ has degree $-1$ if we adopt $\deg u^i=1$.

Consider the linear system 
\begin{align}
&\frac{\p\psi_j}{\p u^i}=\gamma_{ji} \psi_i,\quad i\ne j,\label{w-11}\\
&\frac{\p\psi_i}{\p u^i}=-\sum_{k\ne i} \gamma_{ki} \psi_k.\label{w-12}
\end{align}
The above conditions for $\gamma_{ij}$ ensure the compatibility of this linear system, so its solution space $\S$ has dimention $n$,
and we can find a fundamental system of solutions 
\begin{equation}\label{zh-11-5} 
\Psi_\alpha=(\psi_{1\alpha}(u),\dots, \psi_{n\alpha}(u))^T, \quad  \alpha=1,\dots,n,
\end{equation}
which form a basis of $\S$.

\begin{lem}\label{zh-11-8}
Let $\psi=(\psi_1,\dots, \psi_n)$ be a nontrivial solution of the linear system \eqref{w-11}, \eqref{w-12} on the domain $D$, that is there exist
$i\in \{1, \dots, n\}$, and $u \in D$ such that $\psi_i(u)\ne 0$.
Assume that the rotation coefficients $\gamma_{ij}$
satisfy the irreducibility condition given in Definition \ref{zh-11-2}, then there exists
$u_0\in D$ such that for each $i\in\{1, \dots, n\}$, $\psi_i(u_0)\ne 0$.
\end{lem}
\begin{prf}
For any subset $S\subseteq \{1, \dots, n\}$, define $\phi_S=\prod_{i \in S}\psi_i$.
We assume $\phi_{\{1, \dots, n\}}=0$ on the domain $D$,
then we are to show that $\psi$ is a trivial solution, that is $\phi_{\{i\}}=0$ on $D$ for each
$i=1, \dots, n$. To this end, we will prove that for any $S\subseteq \{1, \dots, n\}$,
$\phi_S=0$  for any $u\in D$ by induction on the size of $S$. We have known that if $\# S=n$,
then $\phi_S=0$. Assume for some $k\le n$, and any $S\subseteq \{1, \dots, n\}$
with $\#S=k$, we have $\phi_S(u)=0$ for any $u\in D$.
For $T\subseteq \{1, \dots, n\}$ with $\# T=k-1$,
and any given $u\in D$, we can find $i\in T$, and $j \notin T$ such that
$\gamma_{ij}(u)\ne 0$ because of the irreducibility condition. Without loss of generality we can assume that $\psi_i(u)$ does not identically vanish. Take $S=T\cup \{j\}$, then consider $\frac{\p \phi_S}{\p u^i}$:
\[0=\frac{\p \phi_S}{\p u^i}=\sum_{k\in S}\phi_{S-\{k\}}\frac{\p \psi_k}{\p u^i}
=\sum_{k\in S, k\ne i}\phi_{S-\{i, k\}}\gamma_{ik}(\psi_i^2-\psi_k^2),\]
so we have
\[\phi_T\frac{\p \phi_S}{\p u^i}=\gamma_{ij}\phi_T^2\psi_i=0.\]
Since $\gamma_{ij}(u)\ne 0$, we have $\phi_T^2\psi_i=0$, which implies
$\phi_T=0$.
\end{prf}

We assume that $\gamma_{ij}$ is irreducible from now on, and shrink $D$ (if necessary) such that $D$ is contractible,
and $\psi_{i1}\ne 0$ on $D$ for each $i=1, \dots, n$.

\begin{lem}We have the following facts:
\begin{itemize}
\item[i)] Define \[\eta_{\alpha\beta}=\sum_{i=1}^n \psi_{i\alpha}\psi_{i\beta},\]
then $(\eta_{\alpha\beta})$ is a constant symmetric non-degenerate matrix. We denote its inverse matrix by $(\eta^{\alpha\beta})$.
\item[ii)] For each $\al=1,\dots, n$, the 1-form
\[\omega_\alpha=\sum_{i=1}^n \psi_{i\alpha} \psi_{i1} d u^i\]
is closed, so there exist smooth functions $v_\alpha$ such that $\omega_\alpha=d v_\alpha$.
Denote $v^\alpha=\eta^{\alpha\beta}v_\beta$, then $(v^1, \dots, v^n)$ can serve as a local coordinate system on $D$. In this local coordinate system we have
\[D_Z=\frac{\p}{\p v^1}.\]
\item[iii)] Define the functions \[c_{\alpha\beta\gamma}=\sum_{i=1}^n \frac{\psi_{i\alpha}\psi_{i\beta}\psi_{i\gamma}}{\psi_{i1}},\]
then $c_{\alpha\beta\gamma}$ are symmetric with respect to the three indices  and satisfy the following conditions:
\begin{align}
&c_{1\alpha\beta}=\eta_{\alpha\beta},\label{cabc-cond-1}\\
& c_{\alpha\beta\xi}\eta^{\xi\zeta} c_{\zeta\gamma\delta}=c_{\delta\beta\xi}\eta^{\xi\zeta} c_{\zeta\gamma\alpha},\label{cabc-cond-2}\\
&\frac{\p c_{\alpha\beta\gamma}}{\p v^\xi}= \frac{\p c_{\xi\beta\gamma}}{\p v^\alpha}. \label{cabc-cond-3}
\end{align}
\end{itemize}
\end{lem}
\begin{prf}
The items i), ii) and the condition \eqref{cabc-cond-1} are easy, so we omit their proofs. The condition \eqref{cabc-cond-2} follows from the identity
$\psi_{i\xi}\eta^{\xi\zeta}\psi_{j\zeta}=\delta_{ij}$. The condition \eqref{cabc-cond-3} can be proved by the chain rule and the following identities
\begin{equation}
\frac{\p v^\alpha}{\p u^i}=\psi_i^\alpha\psi_{i1}, \quad \frac{\p u^i}{\p v^\alpha}=\frac{\psi_{i\alpha}}{\psi_{i1}},
\label{jacobian-uv}
\end{equation}
where $\psi_i^\alpha=\eta^{\alpha\beta}\psi_{i\beta}$.
\end{prf}

The above lemma implies immediately the following corollary.
\begin{cor}\label{cor-potential}
There exists a smooth function $F(v)$ on $D$ such that
\[c_{\alpha\beta\gamma}=\frac{\p^3 F}{\p v^\alpha\p v^\beta\p v^\gamma},\]
and it gives the potential of a Frobenius manifold structure (without the quasi-homogeneity condition) on $D$.
\end{cor}

By using \eqref{jacobian-uv} we have
\[\frac{\p}{\p u^i}\circ \frac{\p}{\p u^j}=c_{\alpha\beta}^\gamma\frac{\p v^\alpha}{\p u^i}\frac{\p v^\beta}{\p u^j}
\frac{\p u^k}{\p v^\gamma}\frac{\p}{\p u^k}=\delta_{ij}\frac{\p}{\p u^i},\]
so $u^1, \dots, u^n$ are the canonical coordinates of this Frobenius manifold. Then its first metric reads
\[\langle d u^i, d u^j\rangle_1=\eta^{\alpha\beta}\frac{\p u^i}{\p v^\alpha}\frac{\p u^j}{\p v^\beta}
=\delta_{ij} \psi_{i1}^{-2},\]
which is in general not equal to the original metric $g_1$ associated to the first 
Hamiltonian structure $P_1$. Though this Frobenius manifold may be not quasi-homogeneous, we can still define its
second metric as follows:
\[\langle d u^i, d u^j\rangle_2=\delta_{ij} u^i \psi_{i1}^{-2}.\]
The two metrics $\langle\ \,,\ \rangle_1$ and $\langle\ \,,\ \rangle_2$ are compatible, since they have the same rotation coefficients with
the original $g_1$, $g_2$ associated to the bihamiltonian structure $(P_1, P_2)$.

The above Frobenius manifold structure depends on the choice of the solution $\Psi_1$ of the linear system \eqref{w-11}, \eqref{w-12}. It is easy to see that
\[\psi_{i1}=f_i^{\frac12}=(f^i)^{-\frac12},\quad i=1,\dots,n\] give a solution to the linear system \eqref{w-11}, \eqref{w-12}.
If we choose it as $\Psi_1$, then the two metrics $\langle\ \,,\ \rangle_1$ and $\langle\ \,,\ \rangle_2$ coincide with $g_1$, $g_2$,
so we call the corresponding Frobenius manifold structure the {\em canonical one} associated to $(P_1, P_2; Z)$.

There are also other choices for $\Psi_1$ such that the corresponding Frobenius manifold is quasi-homogeneous.
By using the identity \eqref{gamma-cond-3}, one can show that Euler vector field $E$ defined by \eqref{euler-vf} acts on the solution space $\S$ as a linear transformation.
Suppose we are working in the complex manifold case, then $E$ has at least one eigenvector in $\S$. We denote this eigenvector by $\Psi_1$,
and denote its eigenvalue by $\mu_1$, then choose other basis $\Psi_2, \dots, \Psi_n$ such that the matrix of $E$ becomes the Jordan normal form,
that is, there exists $\mu_\alpha\in\C$, and $p_\alpha=0$ or $1$, such  that
\[E(\Psi_\alpha)=\mu_\alpha \Psi_\alpha+p_{\alpha-1}\Psi_{\alpha-1}.\]

\begin{lem}
The Frobenius manifold structure corresponding to the above $\Psi_1$ is quasi-homogeneous with the Euler vector field $E$
and the charge $d=-2\mu_1$.
\end{lem}
\begin{prf}
The trivial identity $E(\eta_{\alpha\beta})=0$ implies that
\[\left(\mu_{\alpha}\eta_{\alpha\beta}+p_{\alpha-1}\eta_{(\alpha-1)\beta}\right)
+\left(\mu_{\beta}\eta_{\alpha\beta}+p_{\beta-1}\eta_{\alpha(\beta-1)}\right)=0.\]
Denote by $L_E$ the Lie derivative with respect to $E$, then the identity $L_E \omega_\alpha=d E(v_\alpha)$ implies
\[d E(v_\alpha)=\left(\mu_\alpha+\mu_1+1\right)dv_\alpha+p_{\alpha-1}d v_{\alpha-1},\]
so there exist some constants $r_\alpha\in\C$ such that
\[E(v_\alpha)=\left(\mu_\alpha+\mu_1+1\right)v_\alpha+p_{\alpha-1}v_{\alpha-1}+r_\alpha.\]
On the other hand, we have
\begin{align*}
E(c_{\alpha\beta\gamma})=&\left(\mu_\alpha+\mu_\beta+\mu_\gamma-\mu_1\right)c_{\alpha\beta\gamma}\\
&\quad+p_{\alpha-1}c_{(\alpha-1)\beta\gamma}+p_{\beta-1}c_{\alpha(\beta-1)\gamma}+p_{\gamma-1}c_{\alpha\beta(\gamma-1)}.
\end{align*}
By using the above identities, one can show that
\[\frac{\p^3}{\p v^\alpha\p v^\beta \p v^\gamma}\left(E(F)-(3+2\mu_1) F\right)=0\quad \mbox{for all}\quad \alpha, \, \beta\, \gamma\]
that gives the quasi-homogeneity condition for $F$.
\end{prf}

For each eigenvector $\Psi_1$ of $E$, one can construct a quasi-homogeneous Frobenius manifold. All these Frobenius manifolds (including the
canonical one) are related by Legendre transformations (see \cite{Du-1}).
To see this, let us denote by $F(v)=F(v^1,\dots, v^n)$ and $\tilde{F}(\tilde{v})=\tilde{F}(\tilde{v}^1,\dots, \tilde{v}^n)$ the Frobenius manifold potentials 
constructed above starting from the fundamental solutions 
$(\Psi_1,\dots, \Psi_n)$ and $(\tilde{\Psi}_1,\dots, \tilde{\Psi}_n)$
of the linear system \eqref{w-11}, \eqref{w-12}. These two fundamental solutions are related by a non-degenerate constant matrix $A=(a^\al_\beta)$ by the formula
\[ (\Psi_1,\dots, \Psi_n)=(\tilde\Psi_1,\dots, \tilde\Psi_n) A.\]
Introduce the new coordinates
\[\begin{pmatrix} \hat{v}^1 \\ \vdots \\ \hat{v}^n\end{pmatrix}=A \begin{pmatrix} v^1 \\ \vdots \\ v^n\end{pmatrix}\]
and denote 
\[ \hat{F}(\hat{v})=\hat{F}(\hat{v}^1,\dots, \hat{v}^n):=F(v).\]
Then it is easy to verify that 
\[\hat{v}^\al=\tilde{\eta}^{\al\beta} a^\gamma_1 \frac{\p\tilde{F}(\tilde{v})}{\p \tilde{v}^\beta\p\tilde{v}^\gamma},\quad
\frac{\p^2\hat{F}(\hat{v})}{\p\hat{v}^\al\p\hat{v}^\beta}=\frac{\p^2\tilde{F}(\tilde{v})}{\p\tilde{v}^\al\p\tilde{v}^\beta},\]
and in the $\hat{v}^1,\dots, \hat{v}^n$ coordinates the metrics $g_1, g_2$ 
have the expressions
\[\frac{\p \hat{v}^\al}{\p u^i} g_1^{ij}(u)  \frac{\p \hat{v}^\beta}{\p u^j}
=\tilde{\eta}^{\al\beta}
,\quad
\frac{\p \hat{v}^\al}{\p u^i} g_2^{ij}(u)  \frac{\p \hat{v}^\beta}{\p u^j}
=\tilde{g}(\tilde{v}).\]

\noindent{\em{Proof of Theorem \ref{mainthm00}\,}} 
The first part of the theorem 
follows from the results of \cite{XZ}, and the second part of the theorem is proved by the  arguments given above. The theorem is proved.\hfill{$\Box$}


\section{The principal hierarchy and its tau structure}\label{sec-3}

Let $(P_1, P_2; Z)$ be a flat exact bihamiltonian structure.
Denote 
\begin{align*}
&d_a=\mathrm{ad}_{P_a}:\hF\to\hF,\quad a=1, 2,\\
&\delta=ad_Z:\hF\to\hF.
\end{align*}

\begin{dfn}\mbox{}
\begin{itemize}
\item[i)] Define $\H:=\Ker(d_2\circ d_1)\cap \hF^0$, whose elements are called bihamiltonian conserved quantities.
\item[ii)] Define $\X:=\Ker(d_1)\cap\Ker(d_2)\cap \hF^1$, whose elements are called bihamiltonian vector fields.
\end{itemize}
\end{dfn}
Note that the space $\X$ is actually the bihamiltonian cohomology $BH^1(\hF, P_1, P_2)$, see \cite{BCIH-I}.

\begin{lem}\label{lem-H0X1}
$\H\subset \hF^0_0$, and $\X\subset \hF^1_1$.
\end{lem}
\begin{prf}
If $[P_2, [P_1, H]]=0$, then there exists $K\in \hF^0$ such that $[P_1, H]=[P_2, K]$. By using Lemma 4.1 of \cite{DLZ-1},
we know that $H\in \hF^0_0$.

If $X=\int(X^\alpha\theta_\alpha)\in\hF^1_0$ satisfies $[P_1,X]=[P_2, X]=0$, then we have
\[\nabla_j X^i=0,\quad \tnabla_j X^i=0.\]
Recall that  $\nabla$, $\tnabla$ are the Levi-Civita connections of the metrics $g_1$, $g_2$ associated
with $P_1, P_2$ respectively, 
\begin{equation}\label{zh-12-30a}
\nabla_i=\nabla_{\frac{\p}{\p u^i}},\quad \tnabla_i=\tnabla_{\frac{\p}{\p u^i}}
\end{equation}
and $u^1, \dots, u^n$ are the canonical coordinates of $(P_1, P_2)$.
It follows from the explicit expressions of $g_a^{ij}$, $\Gamma^{ij}_{k,a}$ that $X^i=0$ and so we have
$BH^1_0(\hF, P_1, P_2)\cong 0$.
On the other hand, Lemma 4.1 of \cite{DLZ-1}
implies that $BH^1_{\ge2}(\hF, P_1, P_2)\cong 0$, so consequently $\X=BH^1_1(\hF, P_1, P_2)$.
The lemma is proved.
\end{prf}

\begin{cor}
\begin{itemize}
\item[i)] For any $X, Y \in \X$, we have $[X, Y]=0$;
\item[ii)] For any $X\in \X$, $H\in \H$, we have $[X, H]=0$;
\item[iii)] For any $H, K \in \H$, we have $\{H, K\}_{P_1}:=[[P_1, H], K]=0$.
\end{itemize}
\end{cor}
\begin{prf}
i) If $X, Y \in \X$, then the above lemma shows that $\deg X=\deg Y=1$, so $\deg [X, Y]=2$. But we also have $[X, Y]\in \X$, so $[X, Y]=0$.

ii) If $X\in \X$, $H\in \H$, then $K=[X, H]\in \H$. But $\deg X=1$, $\deg H=0$, so $\deg K=1$, which implies $K=0$.

iii) Take $X=[P_1, H]$, then by applying ii) we obtain $\{H, K\}_{P_1}=0$.
\end{prf}

\begin{lem} \label{lem-23}
We have the following isomorphism
\begin{equation}\label{w-1}
\X\cong \H/\V,
\end{equation}
where $\V=\Ker(d_1)\cap \hF^0$ is the space of Casimirs of $P_1$.
\begin{itemize}
\item[i)] A local functional $H\in\hF^0$ is a bihamiltonian conserved quantity if and only if one can choose its density $h$ so that $h\in \A_0$
and satisfies the condition 
\begin{equation}\label{w-2}
\nabla_i\nabla_jh=0,\quad i\ne j,
\end{equation}
where $\nabla_i=\nabla_{\frac{\p}{\p u^i}}$ are defined as in \eqref{zh-12-30a}.
\item[ii)] A vector field $X\in \hF^1$ is a bihamiltonian vector field if and only if it has the following form
\[X=\int\sum_{i=1}^n A^i(u)u^{i,1}\theta_i,\]
where $A^i(u)$ satisfy the following equations:
\begin{equation}
\frac{\p A^i}{\p u^j}=\Gamma^i_{ij}\left(A^j-A^i\right), \quad \textrm{for } j\ne i, \label{de-for-X}
\end{equation}
here $\Gamma^i_{ij}$ is the Christoffel coefficients of the Levi-Civita connection of $g_1$.
\end{itemize}
\end{lem}
\begin{prf}
Consider the map $\phi=d_1|_{\H}:\H\to \X$. It is easy to see that $\phi$ is well-defined, and $\Ker(\phi)=\V$. Note that 
\[H^1_{\ge 1}(\hF, P_a)\cong 0, \quad a=1,2,\]
so for a given $X\in BH^1_{\ge1}(\hF, P_1, P_2)$, there exists $H, G\in\F$ such that 
\[X=[P_1, H]=[P_2, G].\]
From the second equality we also know that $H\in\H$. So the map $\phi$ is surjective and we proved that the map $\phi$ induces the isomorphism \eqref{w-1}.

Let $H\in\H$, then it yields a bihamiltonian vector field $X=[P_1, H]$. According to Lemma \ref{lem-H0X1}, $H\in\F_0$, $X\in\hF^1_1$.
So we can choose the density of $H=\int(h)$ such that $h\in \A_0$, and
\[X=\int(X^i_ju^{j,1}\theta_i),\]
where $X^i_j=-\nabla^i\nabla_jh$, and $\nabla^i=g^{ik}_1 \nabla_k$. The conditions $[P_1, X]=0$
and  $[P_2, X]=0$ read
\begin{align}
& g_1^{ij}X^k_j=g_1^{kj}X^i_j,\quad \nabla_kX^i_j=\nabla_jX^i_k, \label{X-cond-1}\\
& g_2^{ij}X^k_j=g_2^{kj}X^i_j,\quad \tnabla_kX^i_j=\tnabla_jX^i_k, \label{X-cond-2}
\end{align}
The diagonal form \eqref{zh-10-30} of $g_1$ and $g_2$ and the first equations of \eqref{X-cond-1} and \eqref{X-cond-2} imply that
\[(u^i-u^j) f^j X^i_j=0,\]
so $X^i_j$ is diagonal. Then the second equation of \eqref{X-cond-1} gives the desired equation \eqref{de-for-X}. 
Let $\hat{\Gamma}^i_{ij}$ be the Christoffel coefficients of the Levi-Civita connection of $g_2$, then one can show that for $i\ne j$
\[\hat{\Gamma}^i_{ij}=\Gamma^i_{ij}=\frac{1}{2f_i}\frac{\p f_i}{\p u^j},\]
so the second equation of \eqref{X-cond-2} also gives \eqref{de-for-X}. The lemma is proved.
\end{prf}

\begin{lem}\label{lem-zh-3-3}
We have $\delta(\H)\subseteq \H$. Denote  $\varphi=\delta|_{\H}:\H\to\H$, then $\varphi$ is surjective and $\dim\Ker(\varphi)=n$.
\end{lem}
\begin{prf}
Let $H\in\H$, so we have $[P_2, [P_1, H]]=0$. From the graded Jacobi identity it follows that 
\begin{align*}
[P_2, [P_1, [Z, H]]]=&[P_2, -[[H, P_1],Z]-[[P_1, Z],H]]\\
=& [[P_2,Z],[P_1, H]]+[[P_2, [P_1, H]],Z]=[P_1, [P_1, H]]=0,
\end{align*}
so we have $\delta(\H)\subseteq \H$.

Suppose $H=\int(h)\in\H$,  then from Lemma \ref{lem-23} it follows that the density $h$ can be chosen to belong to $\A_0$ and
$\nabla_i\nabla_j h=0$ for $i\ne j$. If $\varphi(H)=0$, then
\[\sum_{i=1}^n\nabla_i h=0,\]
so we have $\nabla_i\nabla_j h=0$ for any $i, j$, i.e. $h\in\V$. Thus $h$ can be represented as
\[h=\sum_{\alpha=1}^n c_\alpha v^\alpha+c_0,\]
where $c_0, c_1, \dots, c_n$ are some constants, and $v^\alpha$ are the flat coordinates of $g_1$.
From the condition 
$D_Z=\frac{\p}{\p v^1}$ it follows that $c_1=0$, so $\dim\Ker(\varphi)=n$.

To prove that $\varphi$ is surjective, we need to show that for any $g\in\A_0$ satisfying $\nabla_i\nabla_j g=0\ (i\ne j)$, there exists $h\in\A_0$ such that
\begin{equation}\label{zh-08}
\nabla_i\nabla_j h=0\ (i\ne j),\quad \sum_{i=1}^n\nabla_i h=g.
\end{equation}
Denote $\xi_j=\nabla_jh$, then by using the identity \eqref{zh-07}
we know that the above equations imply that
\begin{equation}\label{w-5}
\nabla_i\xi_j=\left\{\begin{array}{cc}0, & i\ne j; \\ \nabla_i g, & i=j.\end{array}\right.
\end{equation}

Let us first prove that the functions $\zeta_{ij}$ defined by the l.h.s. of \eqref{w-5}
satisfy the equalities
\begin{equation}
\nabla_k\zeta_{ij}=\nabla_i\zeta_{kj}. \label{nabla-zeta}
\end{equation}
Denote by $\Gamma_{ij}^k$ the Christoffel coefficients of the first metric,  then we have
\begin{align*}
\nabla_k\zeta_{ij} &= \zeta_{ij,k}-\Gamma^{\alpha}_{ki}\zeta_{\alpha j}-\Gamma^{\alpha}_{kj}\zeta_{i\alpha}\\
&= \zeta_{ij,k}-\Gamma^j_{ki}\zeta_{jj}-\Gamma^i_{kj}\zeta_{ii}.
\end{align*}
Here summation over the repeated upper and lower \emph{Greek} indices is assumed. Note that we do not sum over the repeated \emph{Latin} indices.
Since $\Gamma_{ki}^j=\Gamma_{ik}^j$, in order to prove the identity \eqref{nabla-zeta} we only need to show that
\[\zeta_{ij,k}-\Gamma^i_{kj}\zeta_{ii}=\zeta_{kj,i}-\Gamma^k_{ij}\zeta_{kk}.\]
When $i=j=k$ or $i, j, k$ are distinct, the above equation holds true trivially, so we only need to consider the case when $i=j$ and $i\ne k$.
In this case, the above equation becomes
\[\left(\nabla_i g\right)_{,k}-\Gamma^i_{ki}\nabla_i g+\Gamma^k_{ii}\nabla_k g=0.\]
On the other hand, the function $g$ satisfies $\nabla_k\nabla_i g=0\ (k\ne i)$, which implies
\[\left(\nabla_i g\right)_{,k}=\Gamma^k_{ki}\nabla_k g+\Gamma^i_{ki}\nabla_i g,\]
here we used the fact that $\Gamma_{ij}^k=0$ if $i, j, k$ are distinct. So we only need to show
\[\Gamma^k_{ki}+\Gamma^k_{ii}=0,\quad i\ne k,\] 
which is equivalent to the flatness condition \eqref{jw-5-2}.

The equalities given in  \eqref{nabla-zeta} imply that there exist solutions $\xi_1,\dots, \xi_n$ of the equations \eqref{w-5}. Since $\zeta_{ij}$ are
symmetric with respect to the indices $i, j$, we can find a function $h\in\A_0$ so that
$\xi_i=\nabla_i h$. It follows from \eqref{zh-07} and \eqref{w-5} that $\sum_{i=1}^n \nabla_i h-g$ is a constant, thus by adjusting the function $h$
by adding $c\, v^1$ 
for a certain constant $c$ we prove the existence of $h\in\A_0$ satisfying the equations 
given in \eqref{zh-08}. The lemma is proved.
\end{prf}

The space $\H$ is too big, so we restrict our interest to a ``dense'' (in a certain sense) subspace of $\H$.
\begin{dfn}
Define $\H^{(-1)}=\V$, $\H^{(p)}=\varphi^{-1}\left(\H^{(p-1)}\right)$, and
\[\H^{(\infty)}=\bigcup_{p\ge -1}\H^{(p)}.\]
\end{dfn}

\begin{rmk}
The action of $\varphi$ is just $\frac{\p}{\p v^1}$, so the space $\H^{(\infty)}$ is a polynomial ring in the indeterminate $v^1$.
It is indeed dense in the space of smooth functions in $v^1$ with respect to an appropriate topology.
\end{rmk}

It is easy to see that $\delta(\V)\subseteq \V$, so
\[\V=\H^{(-1)}\subseteq \H^{(0)}\subseteq \cdots\subseteq \H^{(\infty)}.\]
Note that $\dim \H^{(-1)}=n+1$, and
\[\dim \H^{(p)}=\dim\H^{(p-1)}+\dim\Ker (\varphi)=\dim\H^{(p-1)}+n,\]
so we have $\dim \H^{(p)}=n(p+2)+1$.

Suppose the collection of functions 
\[\{h_{\alpha,p}\in\A_0\mid \alpha=1, \dots, n;\ p=0, 1, 2, \dots\}\]
is a calibration of $(P_1, P_2; Z)$ (see Definition \ref{dfn-cali}).
Then it is easy to see that $h_{\alpha, p}\in \H^{(p)}$, and when $p\ge0$, they form a basis of $\H^{(p)}/\H^{(p-1)}$.
When $p=-1$, $\H^{(-1)}=\V$ contains not only $h_{\alpha, 0}=v_{\alpha}$ but also a trivial functional $\int (1)$,
which form a basis of $\H^{(-1)}$. Let us rephrase the conditions that must be satisfied by the functions $h_{\al,p}$ of a calibration as follows:
\begin{flalign}
\qquad & 1.\quad  H_{\alpha,p}=\int(h_{\alpha, p})\in\H,& \label{calib-1} \\
\qquad & 2.\quad h_{\alpha, -1}=v_{\alpha},\quad D_Z(h_{\alpha, p})=h_{\alpha, p-1}\ (p\ge 1),& \label{calib-2} \\
\qquad & 3.\quad \mbox{The normalization condition \eqref{normalize01}}.& \label{calib-3}
\end{flalign}

Now let us proceed to constructing a calibration for the canonical Frobenius manifold structure $F(v)$ of $(P_1, P_2; Z)$.
Following the construction of \cite{Du-1}, we first define the functions
\[\theta_{\alpha,0}(v)=v_\alpha, \quad \theta_{\alpha,1}(v)=\frac{\p F(v)}{\p v^\alpha},\quad  \alpha=1,\dots, n,\]
where $F$ is introduced in Lemma \ref{cor-potential}.
By adding to the function $F(v)$ a certain quadratic term in $v^1,\dots, v^n$, if needed, we can assume that 
\[ \frac{\p^2 F(v)}{\p v^1 \p v^\alpha}=v_\alpha .\]
Thus we have the following relation: 
\[ D_Z\theta_{\alpha,1}=\frac{\p\theta_{\alpha,1}}{\p v^1}=\theta_{\alpha,0}.\]
The functions $\theta_{\alpha,p}(v)$ for $p\ge 2$ can be defined recursively by using the following relations:
\begin{equation}\label{theta-recur}
\frac{\p^2\theta_{\gamma,p+1}(v)}{\p v^\alpha\p v^\beta}=c_{\alpha\beta\xi} \eta^{\xi\zeta} \frac{\p\theta_{\gamma,p}(v)}{\p v^\zeta},\quad \alpha, \beta, \gamma=1,\dots, n.
\end{equation}
The existence of solutions of these recursion relations is ensured by the associativity conditions \eqref{cabc-cond-2}. 
We can require, as it is done in \cite{Du-1},  that these functions also satisfy the following normalization conditions
\[ \frac{\p\theta_\alpha(v; z)}{\p v^\xi} \eta^{\xi\zeta}\frac{\p\theta_{\beta}(v; -z)}{\p v^\zeta}=\eta_{\alpha\beta},\quad
\alpha, \beta=1,\dots, n.\] 
Here $\theta_{\alpha}(v; z)=\sum_{p\ge 0} \theta_{\alpha,p}(v) z^p$. Now we define the functions $h_{\alpha,p}(v)$ so that 
their generating functions $h_\alpha(v; z)=\sum_{p\ge -1}  h_{\alpha,p}(v) z^{p+1}$ satisfy the following defining relation
\begin{equation}
h_{\alpha}(v; z)=\frac1{z} \frac{\p\theta_{\alpha}(v; z)}{\p v^1} -\frac{1}{z} \eta_{\alpha 1}. \label{dfn-h}
\end{equation}
Moreover, these functions also satisfy the normalization condition
\[ \frac{\p h_\alpha(v; z)}{\p v^\xi} \eta^{\xi\zeta}\frac{\p h_{\beta}(v; -z)}{\p v^\zeta}=\eta_{\alpha\beta},\quad
\alpha, \beta=1,\dots, n.\] 
By adding, if needed, a certain linear in $v^1, \dots v^n$ term to the functions $F(v)$ we also have the relations
\begin{equation}\label{zh-09}
h_{\alpha, 0}(v)=\frac{\p F(v)}{\p v^\alpha},\quad \alpha=1,\dots, n.
\end{equation}

For the above constructed functions $\{h_{\alpha, p}\}$, denote $H_{\alpha, p}=\int h_{\alpha, p}$, and define
\begin{equation}\label{w-8}
X_{\alpha, p}=-[P_1, H_{\alpha, p}]=\int \left(\eta^{\gamma\lambda}\p\left(\frac{\delta H_{\alpha, p}}{\delta v^\lambda}\right)\theta_\gamma\right), \quad p\ge 0.
\end{equation} 
Then the associated evolutionary vector field $D_{X_{\alpha, p}}$ (see Definition 2.2 and Equation (2.5) of \cite{BCIH-I} for details)
corresponds to the system of first order quasilinear evolutionary PDEs \eqref{hamilt01}
\begin{equation}\label{jw-5-3}
\frac{\p v^\gamma}{\p t^{\alpha, p}}=D_{X_{\alpha, p}}(v^\gamma),\quad \alpha=1,\dots, n,\,  p\ge 0.
\end{equation}

\begin{lem}\label{lem-cali}
The functions $h_{\alpha,p}$ and the associated local functionals that we constructed above have the following properties:
\begin{itemize}
\item[i)] $H_{\alpha,p}=\int(h_{\alpha, p})\in\H^{(p)}$,
\item[ii)] $h_{\alpha, -1}=v_{\alpha}$, $D_Z(h_{\alpha, p})=h_{\alpha, p-1}\ (p\ge 0)$.
\end{itemize}
\end{lem}
\begin{prf}
According to the definition \eqref{dfn-h} of $h_{\gamma,p}$, we have
\[h_{\gamma, p}=\frac{\p \theta_{\gamma, p+2}}{\p v^1}=\sum_{k=1}^n \frac{\p \theta_{\gamma, p+2}}{\p u^k}.\]
We only need to prove that $H_{\gamma,p}=\int(h_{\gamma, p})\in\H$, that is $\nabla_i\nabla_j h_{\gamma, p}=0$ for $i\ne j$. The other properties are easy to verify.

The condition $\nabla_i\nabla_j h_{\gamma, p}=0$ for $i\ne j$ reads
\[\frac{\p^2 h_{\gamma, p}}{\p u^i\p u^j}=\sum_{l=1}^n\Gamma_{ij}^l\frac{\p h_{\gamma, p}}{\p u^l},\]
which is equivalent to
\begin{equation}
\frac{\p^2 \theta_{\gamma, p+1}}{\p u^i\p u^j}=\sum_{k, l=1}^n\Gamma_{ij}^l\frac{\p^2 \theta_{\gamma, p+2}}{\p u^k \p u^l}.\label{iden-to-prove}
\end{equation}
The recursion relation \eqref{theta-recur} of $\theta_{\alpha, p}$ has the following form in the canonical coordinates:
\[\frac{\p^2 \theta_{\gamma, p+1}}{\p u^i \p u^j}=\delta_{ij}\frac{\p \theta_{\gamma, p}}{\p u^i}+\frac{\p (\psi_i^\alpha\psi_{i1})}{\p u^j}\frac{\p \theta_{\gamma, p+1}}{\p v^\alpha}.\]
Note that $i\ne j$ in the identity \eqref{iden-to-prove}, so its left hand side reads
\[\frac{\p (\psi_i^\alpha\psi_{i1})}{\p u^j}\frac{\p \theta_{\gamma, p+1}}{\p v^\alpha}
=\gamma_{ij}\left(\psi_{i1}\psi_j^\alpha+\psi_{j1}\psi_i^\alpha\right)\frac{\p \theta_{\gamma, p+1}}{\p v^\alpha}.\]
The right hand side of  \eqref{iden-to-prove} then reads
\begin{equation}
\sum_{k, l=1}^n\Gamma_{ij}^l\frac{\p^2 \theta_{\gamma, p+2}}{\p u^k \p u^l}
=\sum_{k, l=1}^n\Gamma_{ij}^l\left(\delta_{kl}\frac{\p \theta_{\gamma, p+2}}{\p u^k}
+\frac{\p (\psi_k^\alpha\psi_{k1})}{\p u^l}\frac{\p \theta_{\gamma, p+2}}{\p v^\alpha}\right). \label{iden-left-hand}
\end{equation}
Note that
\[\sum_{k=1}^n \psi_k^\alpha\psi_{k1}=\delta^\alpha_1\]
is a constant, so the second summation in \eqref{iden-left-hand} vanishes. In the first summation, we have
\[\Gamma_{ij}^l=\gamma_{ij}\left(\delta_{il}\frac{\psi_{j1}}{\psi_{i1}}+\delta_{jl}\frac{\psi_{i1}}{\psi_{j1}}\right), \quad \mbox{for } i\ne j,\]
and
\[\frac{\p \theta_{\gamma, p+2}}{\p u^k}=\frac{\p v^\alpha}{\p u^k}\frac{\p \theta_{\gamma, p+2}}{\p v^\alpha}
=\psi^\alpha_k\psi_{k1}\frac{\p \theta_{\gamma, p+2}}{\p v^\alpha},\]
which leads to the identity \eqref{iden-to-prove}. The lemma is proved.
\end{prf}

\begin{lem} \label{lem-trans}
The first flow $\frac{\p}{\p t^{1,0}}$ is given by the 
translation along the spatial variable $x$, i.e.
\[\frac{\p}{\p t^{1,0}}=\p.\]
\end{lem}
\begin{prf}From our definition \eqref{w-8}, \eqref{jw-5-3} of the evolutionary vector fields we have  
\begin{align*}
&\frac{\p v^\al}{\p t^{1,0}}=\eta^{\al\beta}\p\frac{\p h_{1,0}}{\p v^\beta}=\eta^{\al\beta}\p\frac{\p^2 \theta_{1,2}}{\p v^\beta \p v^1}\\
=&\eta^{\al\beta}\p\frac{\p \theta_{1,1}}{\p v^\beta}=\eta^{\al\beta}\frac{\p^2 \theta_{1,1}}{\p v^\beta \p v^{\gamma}}v^{\gamma}_x=v^{\al}_x.
\end{align*}
Here we use the recursion relation \eqref{theta-recur}. The lemma is proved.
\end{prf}

From Lemma \ref{lem-cali} and Lemma \ref{lem-trans} we have the following proposition.
\begin{prp}\label{whatever}
The collection of functions 
\[\{h_{\alpha, p}(v)\, |\, \alpha=1,\dots, n; p=0, 1,2, \dots\}\]
that we constructed above is a calibration of the flat exact bihamiltonian structure  $(P_1, P_2; Z)$. 
\end{prp}

In the next section, we will use some results proved in the Appendix, which requires that there exists a bihamiltonian vector field
\[X=\int \left(\sum_{i=1}^n A^i(u)u^{i,1}\theta_i\right)\in\X\]
such that for all $i=1, \dots, n$, and for some $u \in D$,
\[\frac{\p}{\p u^i}A^i(u)\ne 0.\]
In this case, $X$ is called \emph{nondegenerate}.

\begin{lem}
If the bihamiltonian vector field $X$ is nondegenerate, then $A^i(u)\ne A^j(u)$ for all $i\ne j$ and for some $u\in D$.
\end{lem}
\begin{prf}
According to \eqref{de-for-X}, if $A^i(u)=A^j(u)$ for $i\ne j$ and $u\in D$, then
\[\frac{\p A^i}{\p u^i}=\frac{\p A^j}{\p u^i}=\Gamma^j_{ji}\left(A^i-A^j\right)=0.\]
The lemma is proved.
\end{prf}

By shrinking the domain $D$, the nondegeneracy condition for $X$ and the result of the above lemma can be modified
to ``for all $u\in D$'' instead of ``for some $u\in D$''.

\begin{lem}\label{new-lemma}
\mbox{}
\begin{itemize}
\item[i)] When $n=1$, the bihamiltonian vector fields $X_{1,p}\ (p>0)$ are always nondegenerate.
\item[ii)] When $n\ge2$, suppose the bihamiltonian structure $(P_1,P_2)$ is irreducible, then there exists a nondegenerate bihamiltonian vector field $X$ satisfying $[Z, X]=0$.
\end{itemize}
\end{lem}
\begin{prf}
We rewrite the bihamiltonian vector field $X_{\alpha, p}$ defined by \eqref{w-8} in the form
\[X_{\alpha, p}=\int\left(\sum_{i=1}^n A^i_{\alpha, p}(u) u^{i,1}\theta_i\right),\]
then $A^i_{\alpha, p}$ satisfy the following equations:
\begin{align*}
&\frac{\p A^i_{\alpha, p}}{\p u^j}=\Gamma^i_{ij}\left(A^j_{\alpha, p}-A^i_{\alpha, p}\right), \quad \mbox{for } j\ne i, \\
&\frac{\p A^i_{\alpha, p}}{\p u^i}=-\sum_{j\ne i}\frac{\p A^i_{\alpha, p}}{\p u^j}+A^i_{\alpha, p-1}, \quad A^i_{1, 0}=1. 
\end{align*}
When $n=1$, we have $A^i_{1, p}=\frac{(u^1)^p}{p!}$, so $X_{1,p}\ (p>0)$ are always nondegenerate.

When $n\ge2$, a bihamiltonian vector field $X=\int\left(\sum_{i=1}^n A^i(u) u^{i,1}\theta_i\right)$ satisfying $[Z, X]=0$ is characterized by the following
equation 
\begin{align}
&\frac{\p A^i}{\p u^j}=\Gamma^i_{ij}\left(A^j-A^i\right), \quad \mbox{for } j\ne i, \label{eq-AA-1}\\
&\frac{\p A^i}{\p u^i}=-\sum_{j\ne i}\frac{\p A^i}{\p u^j}. \label{eq-AA-2}
\end{align}
The solution space of this system has dimention $n$. If $X$ is degenerate, that is, there exists $i_0\in \{1, \dots, n\}$ such that
\[0\equiv\frac{\p A^{i_0}}{\p u^{i_0}}=-\sum_{j\ne i_0}\Gamma^{i_0}_{i_0j}\left(A^j-A^{i_0}\right)=-\sum_{j=1}^n\Gamma^{i_0}_{i_0j}A^j.\]
Since $(P_1, P_2)$ is irreducible, there exists $j_0\in \{1, \dots, n\}$ with $j_0\ne i_0$ such that $\Gamma^{i_0}_{i_0j_0}(u)\ne 0$ for some $u\in D$, so from the above
equation we have
\[A^{j_0}=-\frac{1}{\Gamma^{i_0}_{i_0j_0}}\sum_{k\ne j_0}\Gamma^{i_0}_{i_0k}A^k.\]
Substituting this expression of $A^{j_0}$ into \eqref{eq-AA-1} and \eqref{eq-AA-2}, we obtain a new linear homogeneous system with unknowns $A^k\ (k\ne j_0)$. The dimension of the solution space
of this new system is at most $n-1$, so not all solutions of \eqref{eq-AA-1} and \eqref{eq-AA-2} are degenerate. The lemma is proved.
\end{prf}

\vskip 1em

Let us proceed to prove Proposition \ref{prp-25} which
shows that the functions $h_{\al,p}$ of a calibration of $(P_1, P_2; Z)$ 
satisfy the tau symmetry condition, and the associated principal hierarchy \eqref{hamilt01} possesses Galilean symmetry.

\vskip 1em

\begin{prfn}{Proposition \ref{prp-25}}
By using the chain rule and the properties of $\{h_{\alpha, p}\}$, we have
\begin{align*}
& \frac{\p h_{\alpha, p-1}}{\p t^{\beta, q}}=\sum_{i=1}^n \frac{\p u^i}{\p t^{\beta,q}}\frac{\p h_{\alpha, p-1}}{\p u^i}\\
=& \sum_{i=1}^n \left(\sum_{j=1}^n \nabla^i \nabla_j h_{\beta, q} u^{j,1}\right)
\left(\sum_{k=1}^n \frac{\p^2 h_{\alpha, p}}{\p u^i \p u^k}\right)\\
=& \sum_{i=1}^n \left(f^i \nabla_i \nabla_i h_{\beta, q} u^{i,1}\right)
\left(\sum_{k=1}^n \frac{\p^2 h_{\alpha, p}}{\p u^i \p u^k}\right).
\end{align*}
Note that the flatness of $Z$ implies the identity \eqref{zh-07},
so we have
\[\sum_{k=1}^n \frac{\p^2 h_{\alpha, p}}{\p u^i \p u^k}=\sum_{k=1}^n \nabla_i\nabla_k h_{\alpha, p}
=\nabla_i\nabla_i h_{\alpha, p}.\]
Therefore,
\[\frac{\p h_{\alpha, p-1}}{\p t^{\beta, q}}=\sum_{i=1}^n f^i\left(\nabla_i \nabla_i h_{\beta, q}\right)
\left(\nabla_i\nabla_i h_{\alpha, p}\right) u^{i,1}=\frac{\p h_{\beta, q-1}}{\p t^{\alpha, p}}.\]

Next we show that for any $\gamma=1, \dots, n$,
\[\frac{\p}{\p t^{\alpha, p}}\frac{\p v^\gamma}{\p s}=\frac{\p}{\p s}\frac{\p v^\gamma}{\p t^{\alpha, p}}.\]
The left hand side reads
\[\frac{\p}{\p t^{\alpha, p}}\frac{\p v^\gamma}{\p s}=\frac{\p v^\gamma}{\p t^{\alpha, p-1}}
+\sum_{\beta, q}t^{\beta, q}\frac{\p^2 v^\gamma}{\p t^{\alpha, p}\p t^{\beta, q-1}}.\]
Note that $\frac{\p v^\gamma}{\p t^{\alpha, p}}$ only depends on $v^\mu$ and $v^\mu_x$, so we have
\begin{align*}
&\frac{\p}{\p s}\frac{\p v^\gamma}{\p t^{\alpha, p}}=\frac{\p v^\mu}{\p s}\frac{\p}{\p v^\mu}\left(\frac{\p v^\gamma}{\p t^{\alpha, p}}\right)
+\frac{\p v^\mu_x}{\p s}\frac{\p}{\p v^\mu_x}\left(\frac{\p v^\gamma}{\p t^{\alpha, p}}\right)\\
=&\left(\delta^\mu_1+\sum_{\beta, q}t^{\beta, q}\frac{\p v^\mu}{\p t^{\beta, q-1}}\right)\frac{\p}{\p v^\mu}
\left(\frac{\p v^\gamma}{\p t^{\alpha, p}}\right)\\
&\quad+\left(\sum_{\beta, q}t^{\beta, q}\frac{\p v^\mu_x}{\p t^{\beta, q-1}}\right)\frac{\p v^\mu_x}{\p s}
\frac{\p}{\p v^\mu_x}\left(\frac{\p v^\gamma}{\p t^{\alpha, p}}\right)\\
=&\frac{\p}{\p v^1}\frac{\p v^\gamma}{\p t^{\alpha, p}}+\sum_{\beta, q}t^{\beta, q}\frac{\p^2 v^\gamma}{\p t^{\beta, q-1}\p t^{\alpha, p}},
\end{align*}
so we only need to show that $\frac{\p v^\gamma}{\p t^{\alpha, p-1}}=\frac{\p}{\p v^1}\frac{\p v^\gamma}{\p t^{\alpha, p}}$, which can be easily obtained
from the fact that $X_{\alpha, p-1}=[Z, X_{\alpha, p}]$. The proposition is proved.
\end{prfn}

Since we have $[X_{\beta,q}, H_{\alpha,p-1}]=0$, $\frac{\p h_{\alpha, p-1}}{\p t^{\beta, q}}$ must be a total
$x$-derivative, so there exists a function $\Omega_{\alpha,p;\beta,q}\in \A_0$
such that
\begin{equation}
\frac{\p h_{\alpha, p-1}}{\p t^{\beta, q}}=\frac{\p h_{\beta, q-1}}{\p t^{\alpha, p}}=\p \Omega_{\alpha,p;\beta,q}.
\end{equation}
The functions $\Omega_{\alpha,p;\beta,q}$ are determined up to the addition of constants, so one can adjust the constants such that these functions 
satisfy some other properties which we describe below.

\begin{dfn}\label{zh-12-2}
A collection of functions 
\[\{\Omega_{\alpha,p;\beta,q}\in\A_0\mid \alpha, \beta=1, \dots, n;\ p,q=0, 1, 2, \dots\}\]
is called a tau structure of the flat exact bihamiltonian structure $(P_1, P_2; Z)$
with a fixed calibration $\{h_{\alpha,p}\}$ if the following conditions are satisfied:
\begin{itemize}
\item[i)] $\p\Omega_{\alpha, p; \beta, q}=\frac{\p h_{\alpha, p-1}}{\p t^{\beta, q}}=\frac{\p h_{\beta, q-1}}{\p t^{\alpha, p}}$.
\item[ii)] $\Omega_{\alpha,p;\beta,q}=\Omega_{\beta,q;\alpha,p}$.
\item[iii)] $\Omega_{\alpha,p;1,0}=h_{\alpha,p-1}$.
\end{itemize}
\end{dfn}

\begin{lem}\label{omega-2}
A tau structure $\{\Omega_{\alpha,p;\beta,q}\}$ satisfies the following equations:
\begin{equation}\label{zh-12}
\frac{\p \Omega_{\alpha,p;\beta,q}}{\p t^{\gamma,r}}=\frac{\p \Omega_{\alpha,p;\gamma,r}}{\p t^{\beta,q}}.
\end{equation}
\end{lem}
\begin{prf}By using Definition \ref{zh-12-2} of tau structures we have
\[\p\left(\frac{\p \Omega_{\alpha,p;\beta,q}}{\p t^{\gamma,r}}-\frac{\p \Omega_{\alpha,p;\gamma,r}}{\p t^{\beta,q}}\right)
=\frac{\p}{\p t^{\gamma,r}}\frac{\p h_{\alpha,p-1}}{\p t^{\beta,q}}-\frac{\p}{\p t^{\beta,q}}\frac{\p h_{\alpha,p-1}}{\p t^{\gamma,r}}=0,\]
so the difference between the left hand side and the right hand side of \eqref{zh-12} is a constant. However, both sides can be represented as differential
polynomials of degree $1$, so the constant must be zero. The lemma is proved.
\end{prf}

\begin{dfn}[cf. \cite{DZ-NF}]\label{zh-01-22f}
Let $\{\Omega_{\alpha,p;\beta,q}\}$ be a tau structure of $(P_1, P_2; Z)$ with the calibration $\{h_{\alpha,p}\}$. The family of partial differential equations
\begin{align}
\frac{\p f}{\p t^{\alpha,p}} & =f_{\alpha,p}, \\
\frac{\p f_{\beta,q}}{\p t^{\alpha,p}} & =\Omega_{\alpha,p;\beta,q}(v), \\
\frac{\p v^\gamma}{\p t^{\alpha,p}} & =\eta^{\gamma\xi}\p \Omega_{\alpha,p;\xi,0}(v) 
\end{align}
with unknown functions $(f, \{f_{\beta,q}\}, \{v^\gamma\})$
    is called the tau cover of the principal hierarchy \eqref{hamilt01} with respect to the tau structure $\{\Omega_{\alpha,p;\beta,q}\}$, and the function
$\tau=e^f$ is called the tau function of the principal hierarchy. Here $\al, \beta, \gamma=1,\dots, n, p\ge 0$.
\end{dfn}

By using Lemma \ref{omega-2}, one can easily show that members of the tau cover commute with each other.
It is obvious that the covering map 
\[(f, \{f_{\beta,q}\}, \{v^\gamma\})\mapsto (\{v^\gamma\})\]
pushes forward the
tau cover to the principal hierarchy. This is the reason why it is named ``tau cover''.

In the remaining part of this section, we assume that the calibration $\{h_{\alpha,p}\}$ is constructed from $\{\theta_{\alpha,p}\}$ as above, see
Proposition \ref{whatever}. We can construct, following \cite{Du-1}, the functions $\Omega_{\alpha,p;\beta}(v)$ by
\begin{equation}\frac{\p h_{\alpha}(v; z_1)}{\p v^\xi} \eta^{\xi\zeta} \frac{\p h_{\beta}(v; z_2)}{\p v^\zeta}-\eta_{\alpha\beta}
=(z_1+z_2)\sum_{p, q\ge 0} \Omega_{\alpha, p;\beta, q}(v) z_1^p z_2^q. \label{omega-dfn}
\end{equation}
We can easily prove the following proposition.
\begin{prp}\label{prop-3-12}
The collection of functions 
\[\{\Omega_{\alpha, p; \beta, q}(v)\,|\, \alpha, \beta=1,\dots,; p, q=0,1,2,\dots\}\]
is a tau structure of the exact bihamiltonian structure $(P_1, P_2; Z)$ with the given calibration $\{h_{\al,p}\}$. 
\end{prp}

\begin{lem}\label{omega-3}
The functions $\{\Omega_{\alpha,p;\beta,q}\}$ constructed in \eqref{omega-dfn} satisfy
the identities
\begin{equation}
\frac{\p \Omega_{\alpha,p;\beta,q}}{\p v^1}=\Omega_{\alpha,p-1;\beta,q}+\Omega_{\alpha,p;\beta,q-1}+\eta_{\alpha\beta}\delta_{p0}\delta_{q0}.
\label{omega-iden-1}
\end{equation}
\end{lem}
\begin{prf}
For a fixed pair of indices $\{\al, \beta\}$, the above identities are equivalent to the identity
\begin{equation}
\frac{\p \Omega_{\alpha;\beta}(v; z_1, z_2)}{\p v^1}=(z_1+z_2)\Omega_{\alpha;\beta}(v; z_1, z_2)+\eta_{\alpha\beta}\label{omega-iden-2}
\end{equation}
for the generating function
\[\Omega_{\alpha;\beta}(v; z_1, z_2)=\sum \Omega_{\alpha,p;\beta,q}(v) z_1^p z_2^q.\]
Note that the generation function $h_{\alpha}(v;\,z)$ satisfies
\[\frac{\p h_{\alpha}(v; z)}{\p v^1}=z\,h_{\alpha}(v; z)+\eta_{\alpha 1},\]
then the identity \eqref{omega-iden-2} can be easily proved by using the definition \eqref{omega-dfn}.
The lemma is proved.\end{prf}

\begin{thm}\label{zh-01-22-g}
The tau cover admits the following Galilean symmetry:
\begin{align}
\frac{\p f}{\p s} & =\frac{1}{2}\eta_{\alpha\beta}t^{\alpha,0}t^{\beta,0}+\sum_{\alpha,p}t^{\alpha,p+1}f_{\alpha,p}, \label{SE-1}\\
\frac{\p f_{\beta,q}}{\p s} & =\eta_{\alpha\beta}t^{\alpha,0}\delta_{q0}+f_{\beta,q-1}+\sum_{\alpha,p}t^{\alpha,p+1}\Omega_{\alpha,p;\beta,q}, \label{SE-2}\\
\frac{\p v^\gamma}{\p s} & =\delta^{\gamma}_1+\sum_{\alpha,p}t^{\alpha,p+1}\frac{\p v^{\gamma}}{\p t^{\alpha,p}}. \label{SE-3}
\end{align}
\end{thm}
\begin{prf}
To prove $\frac{\p}{\p s}$ is a symmetry of the tau cover, we only need to show:
\begin{align}
\left[\frac{\p}{\p s}, \frac{\p}{\p t^{\alpha, p}}\right] K=0,\label{identity-to-prove}
\end{align}
where $K=f$, $f_{\beta,q}$, or $v^\gamma$.
Denote the right hand side of \eqref{SE-1} by $W$, then \eqref{SE-2}, \eqref{SE-3} can be written as
\[\frac{\p f_{\beta,q}}{\p s}=\frac{\p W}{\p t^{\beta,q}}, \quad \frac{\p v^\gamma}{\p s}
=\eta^{\gamma\beta}\frac{\p^2 W}{\p t^{1,0}\p t^{\beta,0}},\]
so the identity \eqref{identity-to-prove} is equivalent to the following one:
\[\frac{\p}{\p s}\Omega_{\alpha,p;\beta,q}=\frac{\p^2}{\p t^{\alpha,p}\p t^{\beta,q}}W.\]

By using the chain rule, we have
\begin{align*}
& \frac{\p}{\p s}\Omega_{\alpha,p;\beta,q}=\frac{\p \Omega_{\alpha,p;\beta,q}}{\p v^{\gamma}}\frac{\p v^{\gamma}}{\p s}
=\frac{\p \Omega_{\alpha,p;\beta,q}}{\p v^{\gamma}}
\left(\delta^{\gamma}_1+\sum_{\xi,s}t^{\xi,s+1}\frac{\p v^{\gamma}}{\p t^{\xi,s}}\right)\\
=& \frac{\p \Omega_{\alpha,p;\beta,q}}{\p v^1}+\sum_{\xi,s}t^{\xi,s+1}\frac{\p \Omega_{\alpha,p;\beta,q}}{\p t^{\xi,s}}.
\end{align*}
On the other hand,
\begin{align*}
& \frac{\p^2}{\p t^{\alpha,p}\p t^{\beta,q}}W=\frac{\p}{\p t^{\alpha,p}}\left(
\eta_{\xi\beta}t^{\xi,0}\delta_{q0}+f_{\beta,q-1}+\sum_{\xi,s}t^{\xi,s+1}\Omega_{\xi,s;\beta,q}
\right)\\
=& \eta_{\alpha\beta}\delta_{p0}\delta_{q0}+\Omega_{\alpha,p-1;\beta,q}+\Omega_{\alpha,p;\beta,q-1}
+\sum_{\xi,s}t^{\xi,s+1}\frac{\p \Omega_{\xi,s;\beta,q}}{\p t^{\alpha,p}}.
\end{align*}
The theorem then follows from Lemma \ref{omega-2} and \ref{omega-3}.
\end{prf}

\section{Tau-symmetric integrable Hamiltonian deformations of the principal hierarchy}\label{sec-4}

Let $(P_1, P_2; Z)$ be a flat exact semisimple bihamiltonian structure of hydrodynamic type.
In this and the next section we consider properties of deformations of
the principal hierarchy \eqref{hamilt01} and its tau structure.
To this end, we fix a calibration $\{h_{\alpha, p}\}$ and a tau structure $\{\Omega_{\alpha, p; \beta, q}\}$
as in the previous section, and we assume that $(P_1, P_2; Z)$ is also irreducible.

Note that the principal hierarchy is determined by the first Hamiltonian structure $P_1$ and the calibration $\{h_{\alpha,p}\}$,
so we first consider their deformations.
\begin{dfn}\label{dfn-tau-sym}
The pair $(\tilde{P_1}, \{\tilde{h}_{\alpha, p}\})$ is called a tau-symmetric integrable deformation, or simply a deformation for short,
of $(P_1, \{h_{\alpha, p}\})$ if it satisfies the following conditions:
\begin{itemize}
\item[i)] $\tilde{P}_1\in \hF^2$ has the form
\[\tilde{P}_1=P_1+P^{[2]}_1+P^{[3]}_1+\dots,\]
where $P^{[k]}_1\in \hF^2_{k+1}$, and it is a Hamiltonian structure.
\item[ii)] $\tilde{h}_{\alpha, p}$ has the form
\[\tilde{h}_{\alpha, p}=h_{\alpha, p}+h_{\alpha, p}^{[2]}+h_{\alpha, p}^{[3]}+\cdots,\]
where $h_{\alpha, p}^{[k]}\in\A_{k}$. Define $\tilde{H}_{\alpha, p}=\int (h_{\alpha, p})$, then
for any pair of indices $(\alpha, p)$, $(\beta, q)$ we must have
\begin{equation}
\{\tilde{H}_{\alpha, p}, \tilde{H}_{\beta, q}\}_{\tilde{P}_1}=0.
\end{equation}
Here $\{F, G\}_{\tilde{P}_1}=[[\tilde{P}_1, F], G]$ for $F, G \in \F$.
\item[iii)] Define $\tilde{X}_{\alpha, p}=-[\tilde{P}_1, \tilde{H}_{\alpha, p}]$, and denote $\tilde{\p}_{\alpha, p}=D_{\tilde{X}_{\alpha, p}}$,
then $\{\tilde{h}_{\alpha, p}\}$ satisfy the tau-symmetry condition
\begin{equation}
\tilde{\p}_{\alpha, p}\left(\tilde{h}_{\beta, q-1}\right)=\tilde{\p}_{\beta, q}\left(\tilde{h}_{\alpha, p-1}\right).\label{zh-n-3}
\end{equation}
\end{itemize}
\end{dfn}
\begin{rmk}
Note that we assume the deformation starts from the second degree, i.e. there is no $P_1^{[1]}$ and $h_{\alpha,p}^{[1]}$ terms.
Without this condition we can also prove the next lemma, and then define the tau cover. We add it to avoid some subtle problems in Theorem
\ref{thm-unq-2} (see Remark \ref{tau-rmk} for more details). Note that for integrable hierarchies that arise in the study of semisimple
cohomological field theories, there are no deformations with odd degrees.
\end{rmk}

A deformation of $(P_1, \{h_{\alpha, p}\})$ yields a {\em tau-symmetric  integrable Hamiltonian
deformation of the principal hierarchy} \eqref{hamilt01} which consists of  the flows
\begin{equation}\label{zh-n-17}
\frac{\p v^\al}{\p t^{\beta,q}}=D_{\tilde{X}_{\beta,q}}(v^\al),\quad 1\le \al, \beta\le n,\, q\ge 0.
\end{equation}
Here the evolutionary vector fields are given by
\[\tilde{X}_{\beta,q}=-[\tilde{P}_1, \tilde{H}_{\beta,q}].\]
From the property ii) of Definition \ref{dfn-tau-sym} we know that these deformed 
evolutionary vector fields are mutually commuting, and so the associated flows which we 
denote by $\tilde{\p}_{\beta,q}$ are also mutually commuting. This is the reason why we call the 
above deformed hierarchy \eqref{zh-n-17} an integrable Hamiltonian deformation of the principal hierarchy.
We will show below that the deformed hierarchy also possesses a tau structure. We note that the notion of
\emph{tau-symmetric integrable Hamiltonian deformation} of the principal hierarchy associated to a Frobenius manifold was introduced in \cite{DLYZ}. In the definition given there the following additional conditions are required:
\begin{enumerate}
\item $\tilde{\p}_{1,0}=\p$.
\item $\tilde{H}_{\al,-1}$ are Casimirs of $\tilde{P}_1$.
\end{enumerate}
These two conditions are consequences of the Definition \ref{dfn-tau-sym}. In fact, since the evolutionary vector field $X$
corresponding to the flow $\tilde{\p}_{1,0}-\p$ is a symmetry of the deformed integrable hierarchy and it belongs to
$\hF^1_{\ge 2}$, by using the existence of a non-degenerate bihamiltonian vector field proved in Lemma \ref{new-lemma} and the property ii) of Corollary \ref{app-cor} we know that $X$ must vanishes. Thus we have 
\begin{equation}\label{zh-n-12}
\tilde{\p}_{1, 0}=\p.
\end{equation}
Similarly, from the fact that $[P_1, H_{\al,-1}]=0$ we know that the vector field 
$X=-[\tilde{P}_1, \tilde{H}_{\al,-1}]\in\hF^1_{\ge 2}$. Since it is a symmetry of the deformed integrable hierarchy
\eqref{zh-n-17} we know that it also vanishes. Thus the second condition also holds true.

\begin{lem}
For any deformation $(\tilde{P_1}, \{\tilde{h}_{\alpha, p}\})$ of $(P_1, h_{\alpha, p}\})$, there exists a unique collection of differential polynomials
$\{\tilde{\Omega}_{\alpha, p; \beta, q}\}$
satisfying the following conditions:
\begin{itemize}
\item[i)] 
$\tilde{\Omega}_{\alpha, p; \beta, q}=\Omega_{\alpha, p; \beta, q}+\Omega_{\alpha, p; \beta, q}^{[2]}+\Omega_{\alpha, p; \beta, q}^{[3]}+\cdots$,
where $\Omega_{\alpha, p; \beta, q}^{[k]}\in\A_{k}$.
\item[ii)] $\p \tilde{\Omega}_{\alpha, p; \beta, q}=\tilde{\p}_{\alpha, p}\left(\tilde{h}_{\beta, q-1}\right)$.
\item[iii)] $\tilde{\Omega}_{\alpha, p; \beta, q}=\tilde{\Omega}_{\beta, q; \alpha, p}$, and\,
 $\tilde{\Omega}_{\alpha, p; 1, 0}=\tilde{h}_{\alpha, p-1}$.
\item[iv)] 
$\tilde{\p}_{\gamma, r}\tilde{\Omega}_{\alpha, p; \beta, q}=\tilde{\p}_{\beta, q}\tilde{\Omega}_{\alpha, p; \gamma, r}$.
\end{itemize}
Here $\alpha, \beta, \gamma=1, \dots, n$ and $p, q, r\ge0$.
This collection of differential polynomials $\{\tilde{\Omega}_{\alpha, p; \beta, q}\}$ is called a \emph{tau structure} of $(\tilde{P_1}, \{\tilde{h}_{\alpha, p}\})$.
\end{lem}
\begin{prf}
According to the definition of $\{\tilde{h}_{\alpha, p}\}$,
\[\int\left(\tilde{\p}_{\alpha, p}\left(\tilde{h}_{\beta, q}\right)\right)
=[\tilde{X}_{\alpha, p}, \tilde{H}_{\beta, q}]=-\{\tilde{H}_{\alpha, p}, \tilde{H}_{\beta, q}\}_{\tilde{P}_1}=0,\]
so there exists $\tilde{\Omega}_{\alpha, p; \beta, q}$ satisfying the conditions i), ii). These conditions determine $\tilde{\Omega}_{\alpha, p; \beta, q}$
up to a constant, which has degree zero. Note that the condition i) fixes  the degree zero part of $\tilde{\Omega}_{\alpha, p; \beta, q}$, so it
is unique. The conditions iii) and iv) can be verified by considering the action of $\p$ on both sides of the equalities, as we did 
in the proof of Lemma \ref{omega-2}.
\end{prf}

\begin{dfn}[\cite{DZ-NF}]
The differential polynomials 
\begin{equation}\label{zh-15}
w^{\alpha}=\eta^{\alpha\beta}\tilde{h}_{\beta,-1}=v^{\alpha}+F^{\alpha}_2+F^{\alpha}_3+\cdots,\quad F^\alpha_k\in\A_{k}
\end{equation}
are called the \emph{normal coordinates} of $(\tilde{P_1}, \{\tilde{h}_{\alpha, p}\})$ and of the deformed principal hierarchy \eqref{zh-n-17}. 
\end{dfn}
The properties of the differential polynomials $\tilde{\Omega}_{\al,p;\beta,q}$ enable us to define the tau cover for $(\tilde{P_1}, \{\tilde{h}_{\alpha, p}\})$
and the deformed principal hierarchy \eqref{zh-n-17}, just as we did for the principal hierarchy given in Definition \ref{zh-01-22f}. From \eqref{zh-15} we know
that we can also represent $v^\al$ in the form
\begin{equation}\label{zh-16}
v^\al=w^\al+\tilde{F}^{\alpha}_2+\tilde{F}^{\alpha}_3+\cdots,
\end{equation}
where $\tilde{F}^\alpha_k$ are differential polynomials of $w^1,\dots, w^n$ of degree
$k$. So the functions $\tilde{\Omega}_{\al,p;\beta, q}(v, v_x,\dots)$ can also be represented as  
differential polynomials in $w^1,\dots, w^n$ by the change of coordinates formulae given in \eqref{zh-16}.
\begin{dfn}[c.f. \cite{DZ-NF}]
The family of partial differential equations
\begin{align}
\frac{\p \tilde{f}}{\p t^{\alpha,p}} & =\tilde{f}_{\alpha,p}, \label{zh-18a}\\
\frac{\p \tilde{f}_{\beta,q}}{\p t^{\alpha,p}} & =\tilde{\Omega}_{\alpha,p;\beta,q}, \\
\frac{\p w^\gamma}{\p t^{\alpha,p}} & =\eta^{\gamma\xi}\p \tilde{\Omega}_{\alpha,p;\xi,0}\label{zh-18b}
\end{align}
with the unknowns functions $(\{w^\al\}, \{\tilde{f}_{\alpha,p}\}, \tilde{f})$
is called the tau cover of the deformed principal hierarchy \eqref{zh-n-17} with respect to
the tau structure $\{\tilde\Omega_{\al,p;\beta, q}\}$, and the function $\tilde{\tau}=e^{\tilde{f}}$ is called the tau function of the deformed
principal hierarchy. 
\end{dfn}

\begin{dfn}\label{equivalent}
Suppose $(\tilde{P}_1, \{\tilde{h}_{\alpha, p}\})$ and $(\hat{P}_1, \{\hat{h}_{\alpha, p}\})$ are two deformations of $(P_1, \{h_{\alpha, p}\})$.
Define $\tilde{H}_{\alpha, p}=\int\left(\tilde{h}_{\alpha, p}\right)$ and $\hat{H}_{\alpha, p}=\int\left(\hat{h}_{\alpha, p}\right)$.
If there exists a Miura transformation $e^{\ad_Y}\ (Y\in \hF^1_{\ge1})$ such that
\[\hat{P}_1=e^{\ad_Y}\left(\tilde{P}_1\right), \quad \hat{H}_{\alpha,  p}=e^{\ad_Y}\left(\tilde{H}_{\alpha, p}\right),\]
then we say that $(\tilde{P}_1, \{\tilde{h}_{\alpha, p}\})$ and $(\hat{P}_1, \{\hat{h}_{\alpha, p}\})$ are equivalent.
\end{dfn}

If $(\tilde{P}_1, \{\tilde{h}_{\alpha, p}\})$ and $(\hat{P}_1, \{\hat{h}_{\alpha, p}\})$ are equivalent, then
\[\hat{X}_{\alpha, p}=-[\hat{P}_1, \hat{H}_{\alpha, p}]=-e^{\ad_Y}\left([\tilde{P}_1, \tilde{H}_{\alpha, p}]\right)
=e^{\ad_Y}\left(\tilde{X}_{\alpha, p}\right),\]
which is equivalent to $\hat{\p}_{\alpha, p}=e^{D_Y}\tilde{\p}_{\alpha, p}e^{-D_Y}$.
The associated deformed principal hierarchy has the form (c.f. \eqref{zh-n-17})
\begin{equation}\label{zh-17}
\frac{\p v^\al}{\p t^{\beta,q}}=D_{\hat{X}_{\beta,q}}(v^\al),\quad 1\le \al, \beta\le n,\, q\ge 0.
\end{equation}
It is obtained from \eqref{zh-n-17} by representing the equations of the hierarchy in terms of the new unkown functions
$\tilde{v}^\al=e^{-D_Y}\left(v^\al\right)$ and re-denoting $\tilde{v}^\al, 
\tilde{v}^\al_x,\dots$ by $v^\al, v^\al_x,\dots$.

\begin{thm} \label{thm-unq-2}
Suppose $(\tilde{P}_1, \{\tilde{h}_{\alpha, p}\})$ and $(\hat{P}_1, \{\hat{h}_{\alpha, p}\})$ are two equivalent deformations related
by a Miura transformation $e^{\ad_Y}$, and they have tau structures $\{\tilde{\Omega}_{\al,p;\beta,q}\}$ and $\{\hat{\Omega}_{\al,p;\beta,q}\}$
respectively. Then there exists a differential polynomial $G$ such that
\begin{align*}
& \hat{h}_{\alpha, p}=e^{D_Y}\left(\tilde{h}_{\alpha, p}\right)+\p\hat{\p}_{\alpha, p} G,\\
& \hat{\Omega}_{\alpha, p; \beta, q}=e^{D_Y}\left(\tilde{\Omega}_{\alpha, p; \beta, q}\right)+\hat{\p}_{\alpha, p}\hat{\p}_{\beta, q}G.
\end{align*}
Moreover, suppose $\{\tilde{f}(t), \{\tilde{f}_{\alpha,p}(t)\}, \{\tilde{w}^\alpha(t)\}\}$ is a solution to the tau cover corresponding to the tau structure
$\{\tilde{\Omega}_{\alpha, p; \beta, q}\}$, then
\[\hat{f}(t)=\tilde{f}(t)+G(t),\quad \hat{f}_{\alpha, p}(t)=\tilde{f}_{\alpha, p}+\frac{\p G(t)}{\p t^{\al,p}},\quad
\hat{w}^\alpha(t)=\tilde{w}^\alpha(t)+\eta^{\alpha\beta}\frac{\p^2 G(t)}{\p x\p t^{\beta, 0}}\]
give a solution $\{\hat{f}(t), \{\hat{f}_{\alpha,p}(t)\}, \{\hat{w}^\alpha(t)\}\}$ to the tau cover corresponding to the tau structure
$\{\hat{\Omega}_{\alpha, p; \beta, q}\}$ and the associated deformed principal hierarchy.
Here $G(t)$ is defined from the differential polynomial $G=G(v, v_x,\dots)$ by
\[ G(t)=\left(e^{-D_Y}G(v, v_x,\dots)\right)|_{v^\al=v^\al(\tilde{w}(t)), \tilde{w}_x(t),\dots)},\]
and $v^\al=v^\al(\tilde{w}, \tilde{w}_x,\dots)$ are defined by the relation $\tilde{w}^\al=\eta^{\al, \gamma} \tilde{h}_{\gamma, 0}(v, v_x, \dots)$
just as we did in \eqref{zh-16}.
\end{thm}
\begin{prf}
The condition $\hat{H}_{\alpha,  p}=e^{\ad_Y}\left(\tilde{H}_{\alpha, p}\right)$ implies that there exists $g_{\alpha, p}\in\A_{\ge1}$ such that
\[\hat{h}_{\alpha, p}=e^{D_Y}\left(\tilde{h}_{\alpha, p}\right)+\p g_{\alpha, p}.\]
The tau-symmetry condition $\hat{\p}_{\alpha, p}\hat{h}_{\beta, q-1}=\hat{\p}_{\beta, q}\hat{h}_{\alpha, p-1}$ for 
$\{\hat{h}_{\al,p}\}$ and the one for $\{\tilde{h}_{\alpha, p}\}$
implies that
\[\p\left(\hat{\p}_{\alpha, p}g_{\beta, q-1}-\hat{\p}_{\beta, q}g_{\alpha, p-1}\right)=0,\]
so we have $\hat{\p}_{\alpha, p}g_{\beta, q-1}=\hat{\p}_{\beta, q}g_{\alpha, p-1}$. In particular, by taking $(\beta, q)=(1, 0)$, we have
\[\hat{\p}_{\alpha, p}g_{1,-1}=\p g_{\alpha, p-1},\]
so $\int\left(g_{1,-1}\right)$ gives a conserved quantity for $\hat{\p}_{\alpha, p}$ with a positive degree. According to Theorem \ref{app-thm},
there exists $G\in\A$ such that 
\begin{equation}\label{zh-n-15}
g_{1,-1}=\p G,
\end{equation}
then we have
\[\p\left(\hat{\p}_{\alpha, p} G-g_{\alpha, p-1}\right)=0,\]
so $g_{\alpha,p-1}=\hat{\p}_{\alpha, p}G$ for $\al=1,\dots, n, p\ge 0$. Thus we have
\begin{align*}
&\p \hat{\Omega}_{\alpha, p; \beta, q}=\hat{\p}_{\alpha, p}\hat{h}_{\beta, q-1}=\hat{\p}_{\alpha, p}
\left(e^{D_Y}\left(\tilde{h}_{\beta,q-1}\right)+\p \hat{\p}_{\beta, q} G\right)\\
=& e^{D_Y}\tilde{\p}_{\alpha, p}\left(\tilde{h}_{\beta,q-1}\right)+\p \hat{\p}_{\alpha, p}\hat{\p}_{\beta, q} G
=\p \left(e^{D_Y}\left(\tilde{\Omega}_{\alpha,p;\beta,q}\right)+\hat{\p}_{\alpha, p}\hat{\p}_{\beta, q}G\right),
\end{align*}
so the difference between $\hat{\Omega}_{\alpha, p; \beta, q}$ and $e^{D_Y}\left(\tilde{\Omega}_{\alpha,p;\beta,q}\right)
+\hat{\p}_{\alpha, p}\hat{\p}_{\beta, q}G$ is a constant. However, they have the same leading terms, so the constant must be zero. 

The remaining assertions of the theorem follow from our definition of the tau covers of the deformed principal hierarchies.
The theorem is proved.
\end{prf}

\begin{rmk}\label{tau-rmk}
If in Definition \ref{dfn-tau-sym} we permit the appearance of first degree deformations, i.e. $P_1^{[1]}$ and $h_{\alpha, p}^{[1]}$, the first identity
of the above theorem should be replaced by
\[\hat{h}_{\alpha, p}=e^{D_Y}\left(\tilde{h}_{\alpha, p}\right)+\hat{\p}_{\alpha, p}\sigma,\]
where $\sigma$ is a conserved density of $\hat{\p}_{\alpha, p}$, and the solutions $\hat{f}=\log\hat{\tau}$ and $\tilde{f}=\log\tilde{\tau}$
of the tau covers of  $\hat{\Omega}_{\al,p;\beta,q}$ and $\tilde{\Omega}_{\al,p;\beta,q}$ satisfy the relation
\[\p\left(\log\hat{\tau}-\log\tilde{\tau}\right)=\sigma.\]
The different tau functions defined in \cite{EF, Mira, Wu} for the Drinfeld--Sokolov hierarchies  have such a  relationship. 
\end{rmk}

Next let us consider the Galilean symmetry of the deformed principal hierarchy.
\begin{dfn}
The triple $(\tilde{P}_1, \{\tilde{h}_{\alpha,p}\}, \tilde{Z})$ is a deformation of $(P_1, \{h_{\alpha, p}\}, Z)$ if
\begin{itemize}
\item[i)] The pair $(\tilde{P}_1, \{\tilde{h}_{\alpha,p}\})$ is a deformation of $(P_1, \{h_{\alpha, p}\})$.
\item[ii)] The vector field $\tilde{Z}$ has the form
\[\tilde{Z}=Z+Z^{[2]}+Z^{[3]}+\dots, \quad Z^{[k]}\in \hF^{1}_k,\]
and satisfies conditions $[\tilde{Z}, \tilde{P}_1]=0$ and
\[D_{\tilde{Z}} \tilde{h}_{\al,-1}=\eta_{\al,1},\quad D_{\tilde{Z}} \tilde{h}_{\al,p}=\tilde{h}_{\al,p-1},\quad \al=1,\dots, n,\,p\ge 0.\]
\end{itemize}
\end{dfn}

\begin{lem}\label{omega-w}
Let $\{\tilde{\Omega}_{\alpha,p;\beta,q}\}$ be a tau structure of $(\tilde{P}_1, \{\tilde{h}_{\alpha,p}\}, \tilde{Z})$, and $w^1,\dots, w^n$
are the normal coordinates. Assume that the identity \eqref{omega-iden-1} holds true, then we have:
\[
\frac{\p \tilde{\Omega}_{\alpha,p;\beta,q}}{\p w^1}=
\tilde{\Omega}_{\alpha,p-1;\beta,q}+\tilde{\Omega}_{\alpha,p;\beta,q-1}+\eta_{\alpha\beta}\delta_{p0}\delta_{q0}.
\]
\end{lem}
\begin{prf}
According to Lemma \ref{omega-3}, we only need to show that
\[\p\frac{\p \tilde{\Omega}_{\alpha,p;\beta,q}}{\p w^1}=
\p\tilde{\Omega}_{\alpha,p-1;\beta,q}+\p\tilde{\Omega}_{\alpha,p;\beta,q-1},\]
that is,
\begin{equation}
\frac{\p}{\p w^1}\left(\tilde{\p}_{\beta, q}\left(\tilde{h}_{\alpha, p-1}\right)\right)
=\tilde{\p}_{\beta, q}\left(\tilde{h}_{\alpha, p-2}\right)+\tilde{\p}_{\beta, q-1}\left(\tilde{h}_{\alpha, p-1}\right). \label{idid}
\end{equation}

We first note that one can replace $\frac{\p}{\p w^1}$ by $D_{\tilde{Z}}$. This is  because
\[D_{\tilde{Z}}=\p^s\left(D_{\tilde{Z}}\left(w^\gamma\right)\right)\frac{\p}{\p w^{\gamma, s}}
=\p^s\left(\delta^\gamma_1\right)\frac{\p}{\p w^{\gamma, s}}=\frac{\p}{\p w^1}.\]
Then the identity \eqref{idid} is equivalent to $[D_{\tilde{Z}}, \tilde{\p}_{\beta, q}]=\tilde{\p}_{\beta, q-1}$, which follows from the identities
$\tilde{\p}_{\beta, q}=-D_{[\tilde{P}_1, \tilde{H}_{\beta, q}]}$, and
\[[D_{\tilde{Z}}, D_{[\tilde{P}_1, \tilde{H}_{\beta, q}]}]=D_{[\tilde{Z}, [\tilde{P}_1, \tilde{H}_{\beta, q}]]}
=D_{[\tilde{P}_1, \tilde{H}_{\beta, q-1}]}.\]
The lemma is proved.
\end{prf}

Similar to Theorem \ref{zh-01-22-g}, we have the following theorem on the Galilean
symmetry of the deformed hierarchy $\{\tilde{\p}_{\alpha, p}\}$.
\begin{thm}\label{stringEquation}
Under the assumption of Lemma \ref{omega-w}, the above defined tau cover \eqref{zh-18a}--\eqref{zh-18b} admits the following Galilean symmetry:
\begin{align}
\frac{\p \tilde{f}}{\p s} & =\frac{1}{2}\eta_{\alpha\beta}t^{\alpha,0}t^{\beta,0}+\sum_{\alpha,p}t^{\alpha,p+1}\tilde{f}_{\alpha,p},\\
\frac{\p \tilde{f}_{\beta,q}}{\p s} & =\eta_{\alpha\beta}t^{\alpha,0}\delta_{q0}+\tilde{f}_{\beta,q-1}+\sum_{\alpha,p}t^{\alpha,p+1}\tilde{\Omega}_{\alpha,p;\beta,q}, \\
\frac{\p w^\gamma}{\p s} & =\delta^{\gamma}_1+\sum_{\alpha,p}t^{\alpha,p+1}\frac{\p w^{\gamma}}{\p t^{\alpha,p}}. 
\end{align}
\end{thm}
\begin{prf}
We can prove the theorem by using the same argument as the one given in the proof of Theorem 
\ref{zh-01-22-g}, and by using Lemma \ref{omega-w}.
\end{prf}

\begin{emp}\label{zh-12-31b}
Let $c=\{c_{g,n}: V^{\otimes n}\to H^*(\overline{\mathcal{M}}_{g, n}, \mathbb{Q})\}$ be a semisimple cohomological
field theory. Its genus zero part defines a semisimple Frobenius manifold, which corresponds to a flat exact semisimple bihamiltonian structure of hydrodynamic type. Its principal hierarchy
has a useful deformation, called topological deformation, such that the partition function of $c$ is a tau function of this deformed hierarchy
\cite{DZ-NF, BPS-1, BPS-2}. On the other hand, Buryak constructed another deformation, called double ramification deformation, from the same data,
and conjectured that they are actually equivalent  \cite{Bu}. This conjecture is refined in \cite{BD}  as follow:

Suppose $\F$ is the free energy of the topological deformation. Buryak \emph{et al} show that there exists a unique differential polynomial $P$
such that $\F^{\mathrm{red}}=\F+P$ satisfies the following condition:
\[\left.\mathrm{Coef}_{\epsilon^{2g}}\frac{\p ^{n}\F^{\mathrm{red}}}{\p t^{\alpha_1, p_1}\cdots\p t^{\alpha_n, p_n}}\right|_{t^{*, *}=0}=0,\quad 
p_1+\cdots+p_n\le 2g-2. \]
It is conjectured that $\F^{\mathrm{red}}$ is just the free energy of the double ramification deformation.

Buryak \emph{et al}'s refined conjecture is compatible with our Theorem \ref{thm-unq-2}. They also show that the double ramification deformation
satisfies the string equation, which can also be derived from our Theorem \ref{stringEquation}.
\end{emp}

\section{Tau-symmetric bihamiltonian deformations of the principal hierarchy}\label{sec-5}

In this section, we construct a class of tau-symmetric integrable Hamiltonian deformations 
of the principal hierarchy associated with a semisimple flat exact bihamiltonian structure $(P_1, P_2; Z)$ of hydrodynamic type.
These deformations of the principal hierarchies are in fact bihamiltonian integrable hierarchies.

From \cite{CPS-2, BCIH-I} we know that the bihamiltonian structure $(P_1, P_2)$ 
possesses deformations of the form
\[ \tilde{P}_1=P_1+\sum_{k\ge 1}  Q_{1,k},\quad 
\tilde{P}_2=P_2+\sum_{k\ge 1}  Q_{2,k},\quad  Q_{1,k}, Q_{2, k}\in \hF^2_{k+1}\]
such that $(\tilde{P}_1, \tilde{P}_2)$ is still a bihamiltonian structure, i.e.
\[[\tilde{P}_a, \tilde{P}_b]=0, \quad a, b=1, 2.\]
The space of deformations of the bihamiltonian structure $(P_1, P_2)$ is characterized
by the central invariants $c_1(u), \dots,  c_n(u)$ of $(\tilde{P}_1, \tilde{P}_2)$.
The following theorem of Falqui and Lorenzoni gives a condition under which the deformed bihamiltonian structure inherits the exactness property.
This means that there exists a vector field
$\tilde{Z}\in \hF^1$ such that 
\[[\tilde{Z}, \tilde{P}_1]=0,\quad [\tilde{Z}, \tilde{P}_2]=\tilde{P}_1.\]

\begin{thm}[\cite{FL}]\label{thm-FL}
The deformed bihamiltonian structure $(\tilde{P}_1, \tilde{P}_2)$ is exact if and only if its central invariants $c_1, \dots, c_n$ are constant functions.
Moreover, there exists a Miura type transformation $g$ such that 
\begin{equation}\label{zh-16-1} 
g(\tilde{P}_1)=P_1,\quad 
g(\tilde{P}_2)=P_2+\sum_{k\ge 1}  Q_{k},\quad  Q_{k}\in \hF^2_{2k+1}
\end{equation}
and  $g(\tilde{Z})=Z$, where $Z=Z_0$ is given by \eqref{eq-Z0}.
\end{thm}

In what follows, we assume that $(\tilde{P}_1, \tilde{P}_2; \tilde{Z})$ is a deformation of the
flat exact bihamiltonian structure $(P_1, P_2; Z)$ with constant central invariants $c_1, \dots, c_n$, 
$\tilde{P}_1, \tilde{P}_2$ have the form given in \eqref{zh-16-1}, and $\tilde{Z}=Z$.
We denote by $u^1,\dots, u^n$ and $v^1,\dots, v^n$ the canonical coordinates of $(P_1, P_2)$ and the flat coordinates of $P_1$ 
respectively. We also fix a calibration 
\[\{h_{\alpha, p}(v)\in\A_0 \mid \alpha=1,\dots, n;\ p=0, 1,2, \dots\}\]
and a tau structure 
\[\{\Omega_{\alpha,p;\beta,q}(v)\in\A_0\mid \alpha, \beta=1, \dots, n;\ p,q=0, 1, 2, \dots\}\]
 of the flat exact bihamiltonian structure $(P_1, P_2; Z)$ (see above their construction given in Propositions \ref{whatever}, \ref{prop-3-12}).

We define the space of Casimirs of $\tilde{P}_1$, the space of bihamiltonian conserved quantities and the space of bihamiltonian vector fields respectively,
just like we did for $(P_1, P_2)$, as follows:
\begin{align*}
\tV:=& \Ker([\tilde{P}_1, \cdot])\cap \F, \\
\tH:=& \Ker([\tilde{P}_2, [\tilde{P}_1, \cdot]])\cap \F, \\
\tX:=& \Ker([\tilde{P}_1, \cdot])\cap \Ker([\tilde{P}_2, \cdot])\cap \hF^1.
\end{align*}

\begin{thm}\label{thm-31}
We have the following isomorphisms:
\begin{equation}
\V\cong \tV,\quad \H\cong \tH,\quad \X\cong \tX. 
\end{equation}
In particular, $\tX\cong\tH/\tV$.
\end{thm}
\begin{prf}
Since $\tilde{P}_1=P_1$, we only need to prove that $\H\cong \tH$, $\X\cong \tX$.
Suppose $H\in\tH$ is a bihamiltonian conserved quantity of $(\tilde{P}_1, \tilde{P}_2)$. Expand $H$ as the sum of homogeneous components
\[H=H_0+H_1+H_2+\cdots,\quad H_k\in\F_{2 k},\]
then $H_0$ is a bihamiltonian conserved quantity of $(P_1, P_2)$, so we have a map $\pi:\tH\to\H$, $H\mapsto H_0$.
The fact that $\H$ is concentrated in degree zero (see Lemma \ref{lem-23}) implies that $\pi$ is injective. To prove the isomorphism $\H\cong\tH$,
we only need to show that $\pi$ is surjective, that is, for any bihamiltonian conserved quantity $H_0$ of $(P_1, P_2)$ there exists a
bihamiltonian conserved quantity $H$ of $(\tilde{P}_1, \tilde{P}_2)$ with $H_0$ as its leading term.

Recall that $(\tilde{P}_1, \tilde{P}_2; Z)$ takes the form \eqref{zh-16-1}. If we 
denote $d_a=[P_a, \cdot]\ (a=1,2)$, then $Q_k$ satisfy the following equations:
\[d_1Q_k=0,\quad d_2 Q_k+\frac12\sum_{i=1}^{k-1}[Q_i, Q_{k-i}]=0.\]
We assert that, for any bihamiltonian conserved quantity $H_0\in\H$ of $(P_1, P_2)$, there exists $H_k\in\F_{2k}$ such that
\[H=H_0+H_1+H_2+\cdots\]
is a bihamiltonian conserved quantity of $(\tilde{P}_1, \tilde{P}_2)$. This assertion is equivalent to the solvability of the following equations
for $H_k$ to be solved recursively:
\[d_1d_2H_k=\sum_{i=1}^k[Q_i, d_1 H_{k-i}], \quad k=1, 2, \dots.\]
Assume that we have already solved the above equations for $H_1, \dots, H_{k-1}$ starting from $H_0$. Denote by $W_k$ the right hand side of the above equation.
Then it is easy to see that $d_1W_k=0$, and
\begin{align*}
 & d_2W_k = \left[P_2, \sum_{i=1}^k[Q_i, d_1 H_{k-i}]\right]\\
=& -\sum_{i=1}^k\left([[d_1H_{k-i}, P_2], Q_i]+[[P_2, Q_i], d_1 H_{k-i}]\right)\\
=& \sum_{i=1}^k\left([d_1d_2H_{k-i}, Q_i]+[-d_2Q_i, d_1 H_{k-i}]\right)\\
=& \sum_{i=1}^k\sum_{j=1}^{k-i}[[Q_j, d_1 H_{k-i-j}], Q_i]+\frac12\sum_{m=1}^k\sum_{i=1}^{m-1}[[Q_i,Q_{m-i}],d_1H_{k-m}]\\
=& \frac12\sum_{i,j\ge1, l\ge0, i+j+l=k}\left([[Q_j, d_1 H_l], Q_i]+[[Q_i, d_1 H_l], Q_j]+[[Q_i,Q_j], d_1 H_l]\right)\\
=& 0,
\end{align*}
so $W_k\in\Ker(d_1)\cap\Ker(d_2)\cap\hF^2_{\ge4}$. Since $BH^2_{\ge4}(\hF)\cong0$, there exists $H_k\in\F$ such that
$W_k=d_1d_2H_k$. Thus the isomorphism $\H\cong\tH$ is proved.

It is easy to see that the map 
\[ d_1:\, \tH/\tV\to\tX,\quad H\mapsto X=-[\tilde{P}_1, H]\]
gives the isomorphism $\tH/\tV\cong \tX$, which also induces the isomorphism $\X\cong\tX$.
The theorem is proved.
\end{prf}

It follows from the above theorem that there exist unique deformations 
\[\tilde{H}_{\al, p}=H_{\al,p}+H_{\al,p}^{[1]}+H_{\al,p}^{[2]}+\dots,\quad
H_{\al,p}^{[k]}\in \F_{2k}\]
of the bihamiltonian conserved quantities $H_{\al, p}=\int (h_{\al,p})\in \H$ such that, together 
with the constant local functional $\int (1)$, they form a basis of the subspace
\[ \tH^{\infty}=\bigcup_{p\ge0}\tH^{(p)}\]
of $\tH$, where $\tH^{(p)}$ is the image of $\H^{(p)}$ in $\tH$ of the isomorphism given in the above theorem. For any pair of indices $(\al, p)$, $(\beta, q)$,
it is easy to see that  the local functional 
$H=\{\tilde{H}_{\al,p}, \tilde{H}_{\beta,q}\}_{\tilde{P}_1}:=[[\tilde{P}_1, \tilde{H}_{\al,p}], \tilde{H}_{\beta,q}]$ is a bihamiltonian conserved quantity
w.r.t. $(\tilde{P}_1, \tilde{P}_2)$. Since $H\in\F_{\ge 1}$ we obtain 
\begin{equation}\label{zh-n-19}
\{\tilde{H}_{\al,p}, \tilde{H}_{\beta,q}\}_{\tilde{P}_1}=0
\end{equation}
by using Lemma \ref{new-lemma} and the property i) of Corollary \ref{app-cor}.

Define an operator
\begin{equation}\label{zh-21}
\delta_Z:\hat\F\to\hat\A, \quad Q\mapsto \sum_{i=1}^n \frac{\delta Q}{\delta u^i}=\frac{\delta{Q}}{\delta v^1}.
\end{equation}
Here we used the fact that 
\[ D_Z=\frac{\p}{\p v^1}=\sum_{i=1}^n \frac{\p}{\p u^i}.\]
Then for a local functional $H\in\F$ we have $[Z, H]=\int\left(\delta_Z(H)\right)$. 
Now let us define 
\begin{equation} \label{var-h}
\tilde{h}_{\al,p}=\delta_Z \tilde{H}_{\al,p+1},\quad \al=1,\dots, n, \ p=-1, 0,1,\dots.
\end{equation}

\begin{thm}\label{zh-01-22-a}
The triple $(\tilde{P}_1, \{\tilde{h}_{\alpha,p}\}, \tilde{Z})$ gives a deformation of $(P_1, \{h_{\alpha, p}\}, Z)$.
\end{thm}
\begin{prf}
Define $\tilde{H}'_{\alpha, p}=\int\left(\tilde{h}_{\alpha, p}\right)$.
From the definition of $\tilde{h}_{\al,p}$ we see that $\tilde{H}'_{\alpha, p}=[Z, \tilde{H}_{\al,p+1}]$, so it belongs to
$\tH$.
From the property $D_Z h_{\al, p+1}=h_{\al,p}$ we know that $\tilde{H}'_{\alpha, p}$ and
$\tilde{H}_{\al,p}$ have the same leading term $\int(h_{\al,p})$. Since the bihamiltonian 
conserved quantities of $(\tilde{P}_1, \tilde{P}_2)$ are uniquely determined by their leading terms,
we obtain 
\[\tilde{H}'_{\alpha, p}=\tilde{H}_{\al,p}.\]
In particular, we know from \eqref{zh-n-19} that $\{\tilde{H}'_{\alpha, p}, \tilde{H}'_{\beta, q}\}_{\tilde{P}_1}=0$, and
\[\tilde{X}'_{\al,p}=-[\tilde{P}_1, \tilde{H}'_{\alpha, p}]=-[\tilde{P}_1, \tilde{H}_{\alpha, p}]=\tilde{X}_{\al,p}.\]

Denote by $\bar{\theta}_\alpha$ the super variables corresponding to the flat coordinates $v^1, \dots, v^n$. Recall that
\[\tilde{P}_1=P_1=\frac{1}{2}\int\left(\eta^{\alpha\beta}\bar{\theta}_\alpha\bar{\theta}_\beta^1\right),
\quad \mbox{where } \eta^{\alpha\beta}=\langle dv^\alpha, dv^\beta\rangle_{g_1},\]
so we have
\[\tilde{X}_{\al,p}=\int\left(\eta^{\beta\gamma}\p\left(\frac{\delta \tilde{H}_{\al,p}}{\delta v^\gamma}\right)\bar{\theta}_\beta\right).\]
Denote by $V=\eta_{1\gamma} v^\gamma$, then
\[\frac{\p V}{\p t^{\al,p}}=\frac{\delta \tilde{X}_{\al,p}}{\delta \bar{\theta}_\gamma}\frac{\p V}{\p v^\gamma}
=\eta_{1\gamma}\eta^{\gamma\beta}\p \left(\frac{\delta \tilde{H}_{\al,p}}{\delta v^\beta}\right)=
\p\left(\frac{\delta \tilde{H}_{\al,p}}{\delta v^1}\right),\]
which implies that
\[\p\left(\tilde{\p}_{\alpha, p}\left(\tilde{h}_{\beta, q-1}\right)-\tilde{\p}_{\beta, q}\left(\tilde{h}_{\alpha, p-1}\right)\right)
=\frac{\p}{\p t^{\al,p}}\frac{\p V}{\p t^{\beta,q}}-\frac{\p}{\p t^{\beta,q}}\frac{\p V}{\p t^{\al,p}}=0.\]
Since the difference $\tilde{\p}_{\alpha, p}\left(\tilde{h}_{\beta, q}\right)-\tilde{\p}_{\beta, q}\left(\tilde{h}_{\alpha, p}\right)$
is a differential polynomial with terms of degree greater or equal to one, so it must be zero. The above computation shows that
$(\tilde{P}_1, \{\tilde{h}_{\alpha, p}\})$ is a deformation of $(P_1, \{h_{\alpha, p}\})$, see Definition \ref{dfn-tau-sym}.

Next let us consider the action of $D_Z$ on $\tilde{h}_{\al,p}$. We have
\[D_Z(\tilde{h}_{\al,p+1})
=\frac{\p}{\p v^1}\frac{\delta}{\delta v^1}\tilde{H}_{\alpha, p+2}
=\frac{\delta}{\delta v^1}\frac{\delta}{\delta v^1}\tilde{H}_{\alpha, p+2}
=\frac{\delta}{\delta v^1}\tilde{H}_{\alpha, p+1}=\tilde{h}_{\alpha, p}.\]
Here we used the following identity for variational derivatives:
\[\frac{\p}{\p v^1}\frac{\delta}{\delta v^1}=\frac{\delta}{\delta v^1}\frac{\delta}{\delta v^1},\]
which is a particular case of the identity (i) of Lemma 2.1.5 in \cite{Jacobi}.

We still need to check the identities $D_Z(\tilde{h}_{\al,-1})=\eta_{\alpha, 1}$, which is equivalent to $\delta_Z \tilde{H}_{\alpha,-1}=\eta_{\alpha, 1}$.
Note that the leading term $\tilde{H}_{\alpha,-1}^{[0]}=\int \left(\eta_{\alpha\beta}v^\beta\right)$ of $\tilde{H}_{\alpha,-1}$ is a Casimir of $P_1=\tilde{P_1}$,
so it also belongs to $\tH$. On the other hand, elements of $\tH$ are determined by their leading terms, so we have
$\tilde{H}_{\alpha,-1}=\tilde{H}_{\alpha,-1}^{[0]}$, which implies the desired identity. The theorem is proved.
\end{prf}

\begin{rmk}
Our construction \eqref{var-h} of the Hamiltonian densities that satisfy the tau
symmetry property follows the approach given in \cite{DZ-NF} for the construction of the tau structure of the KdV hierarchy.
Note that this approach was also employed in \cite{BD} to construct tau structures for the double ramification hierarchies
associated to cohomological field theories.
\end{rmk}

The deformation $(\tilde{P}_1, \{\tilde{h}_{\alpha, p}\}, \tilde{Z})$ constructed in the above theorem depends on the choice of $\tilde{P}_2$.
It is natural to ask: if we start from another deformation $(\hat{P}_1, \hat{P}_2; \hat{Z})$ which has the same central invariants as
$(\tilde{P}_1, \tilde{P}_2; \tilde{Z})$ does, how does the result on the deformation $(\tilde{P}_1, \{\tilde{h}_{\alpha, p}\}, \tilde{Z})$ change?

Without loss of generality, we can assume that both $(\hat{P}_1, \hat{P}_2; \hat{Z})$ and $(\tilde{P}_1, \tilde{P}_2; \tilde{Z})$ have been transformed
to the form \eqref{zh-16-1}. If $(\hat{P}_1, \hat{P}_2)$ has the same central invariants as $(\tilde{P}_1, \tilde{P}_2)$,
then there exists a Miura type transformation of the second type
\[\mathrm{v}\mapsto \bar{\mathrm{v}}=e^{-D_Y}\left(\mathrm{v}\right)\]
with $Y\in\hF^{1}_{\ge 2}$ such that
\[\hat{P}_a=e^{\ad_{Y}}\left(\tilde{P}_a\right), \quad a=1, 2.\]
Note that $\hat{P}_1=\tilde{P}_1=P_1$, so $[P_1, Y]=0$, which implies that there exists $K\in\F_{\ge1}$ such that
$Y=[P_1, K]$.

\begin{lem}
The vector field $Y$ and the functional $K$ satisfy $[Y, Z]=0$ and $[K, Z]=0$.
\end{lem}
\begin{prf}
Denote $Z'=e^{\ad_{Y}}\left(Z\right)$, then we have
\begin{align*}
&[\hat{P}_1, Z']=e^{\ad_{Y}}\left([\tilde{P}_1, Z]\right)=0,\\
&[\hat{P}_2, Z']=e^{\ad_{Y}}\left([\tilde{P}_2, Z]\right)=\hat{P}_1,
\end{align*}
so $W=Z'-Z$ is a bihamiltonian vector field of $(\hat{P}_1, \hat{P}_2)$. On the other hand, $W\in\hF^1_{\ge2}$, so we have $W=0$ and, consequently, we have
$[Y, Z]=0$.

It follows from the identity $[Y,Z]=0$ that $[P_1, [K, Z]]=0$, so $C=[K, Z]$ is a Casimir of $P_1$. Since $C\in\F_{\ge1}$, we obtain $C=0$. The lemma is proved.
\end{prf}

From the above lemma we have
\[\int\left(\delta_Z K\right)=[Z,K]=0,\]
so there exists $g\in\A$ such that 
\begin{equation}\label{zh-20}
\delta_Z K=\p g.
\end{equation}

Let $\{\hat{H}_{\alpha,p}\}$, $\{\tilde{H}_{\alpha,p}\}$ be the bihamiltonian conserved quantities of $(\hat{P}_1, \hat{P}_2)$ and 
$(\tilde{P}_1, \tilde{P}_2)$ respectively with the same leading terms $\{h_{\alpha,p}\}$, and $\{\hat{X}_{\alpha,p}\}$, $\{\tilde{X}_{\alpha,p}\}$ be
the corresponding bihamiltonian vector fields:
\[\hat{X}_{\alpha, p}=-[P_1, \hat{H}_{\alpha,p}], \quad \tilde{X}_{\alpha, p}=-[P_1, \tilde{H}_{\alpha,p}].\]
They are related by
\[\hat{H}_{\alpha,p}=e^{\ad_{Y}}\left(\tilde{H}_{\alpha,p}\right), \quad
\hat{X}_{\alpha,p}=e^{\ad_{Y}}\left(\tilde{X}_{\alpha,p}\right).\]
The flows corresponding to $\{\hat{X}_{\alpha,p}\}$ and $\{\tilde{X}_{\alpha,p}\}$ are denoted respectively by $\{\hat{\p}_{\alpha, p}\}$ and
$\{\tilde{\p}_{\alpha, p}\}$. We also have the associated triples $(\tilde{P}_1, \{\tilde{h}_{\alpha,p}\}, \tilde{Z})$ and
$(\hat{P}_1, \{\hat{h}_{\alpha,p}\}, \hat{Z})$ which are constructed in Theorem \ref{zh-01-22-a}. Let $\{\tilde{\Omega}_{\al,p;\beta,q}\}$ and 
$\{\hat{\Omega}_{\al,p;\beta,q}\}$ be the corresponding tau structures. Then the relation between these tau structures and the solutions of the associated
tau covers of the deformed principal hierarchies is given by Theorem \ref{thm-unq-2}, and the following theorem gives the explicit expression of the
differential polynomial $G$.

\begin{thm}\label{thm-unq-1}
The differential polynomial $G$ of Theorem \ref{thm-unq-2} is given by the formula
\begin{equation}\label{zh-19}
G=\sum_{i=1}^\infty \frac{1}{i!}D_Y^{i-1}\left(g\right),
\end{equation}
where the function $g$ is defined in \eqref{zh-20}.
\end{thm}
\begin{prf}
From our construction of the densities of the Hamiltonians we have
\[\hat{h}_{\alpha,p}=\delta_Z \hat{H}_{\alpha,p+1},\quad \tilde{h}_{\alpha,p}=\delta_Z \tilde{H}_{\alpha,p+1},\]
so
\[\hat{h}_{\alpha,p}=\delta_Z\left(e^{\ad_{Y}}\left(\tilde{H}_{\alpha,p+1}\right)\right)
=\sum_{k=0}^\infty \frac{1}{k!}\delta_Z\left(\ad_{Y}^k\left(\tilde{H}_{\alpha,p+1}\right)\right).\]
By using the the definition \eqref{zh-21} of $\delta_Z$ and the identities 
given in Lemma \ref{last-lem} we can show that
\[\delta_Z\left(\ad_{Y}^k\left(\tilde{H}_{\alpha,p+1}\right)\right)=D_Y^k\left(\tilde{h}_{\alpha,p}\right)
+\sum_{i=1}^k\binom{k}{i}D_{\ad_Y^{k-i}\left(\tilde{H}_{\alpha,p+1}\right)} D_Y^{i-1}\left(\delta_Z Y\right),\]
so we have
\begin{align*}
\hat{h}_{\alpha,p}
=&\sum_{k=0}^\infty \frac{1}{k!}\left(D_Y^k\left(\tilde{h}_{\alpha,p}\right)
+\sum_{i=1}^k\binom{k}{i}D_{\ad_Y^{k-i}\left(\tilde{H}_{\alpha,p+1}\right)} D_Y^{i-1}\left(\delta_Z Y\right)\right)\\
=&e^{D_Y}\left(\tilde{h}_{\alpha,p}\right)+D_{\hat{H}_{\alpha, p+1}}\left(\sum_{i=1}^\infty \frac{1}{i!}D_Y^{i-1}\left(\delta_Z Y\right)\right)
\end{align*}
By using the fact that
\[\delta_Z Y=\delta_Z[P_1, K]=D_{P_1}(\delta_Z K)=\p D_{P_1}(g),\]
and $[D_Y, D_{P_1}]=0$, $[\p, D_{Q}]=0$ for $Q\in\hat\F$ (see Lemma \ref{last-lem}), we obtain
\[\hat{h}_{\alpha,p}=e^{D_Y}\left(\tilde{h}_{\alpha,p}\right)+\p D_{\hat{H}_{\alpha, p+1}}D_{P_1}G,\]
where
\[G=\sum_{i=1}^\infty \frac{1}{i!}D_Y^{i-1}\left(g\right).\]
Then by using the identity (see Lemma \ref{last-lem}) 
\[D_{H}D_{P_1}=D_{-[P_1, H]}-D_{P_1}D_{H},\]
and the fact that $D_{H}\left(G\right)=0$, we have
\[\hat{h}_{\alpha,p}=e^{D_Y}\left(\tilde{h}_{\alpha,p}\right)+\p D_{\hat{X}_{\alpha, p+1}}G
=e^{D_Y}\left(\tilde{h}_{\alpha,p}\right)+\p \hat{\p}_{\alpha, p} G.\]
The theorem is proved.
\end{prf}

In the proof of the above theorem the following lemma is used.
\begin{lem}\label{last-lem}
The operator 
\[D_P=\sum_{s\ge0}\left(\p^s\left(\frac{\delta P}{\delta \theta_\alpha}\right)\frac{\p}{\p u^{\alpha,s}}
+(-1)^p\p^s\left(\frac{\delta P}{\delta u^\alpha}\right)\frac{\p}{\p \theta_\alpha^s}\right),\quad P\in\hF^p\]
and the bracket
\[[P, Q]=\int\left(\frac{\delta P}{\delta \theta_\alpha}\frac{\delta Q}{\delta u^\alpha}
+(-1)^p\frac{\delta P}{\delta u^\alpha}\frac{\delta Q}{\delta \theta_\alpha}\right), \quad P\in \hF^p,\ Q\in\hF^q\]
satisfy the following identities:
\begin{align}
&  [\p, D_P]=0;\nn\\
& \frac{\delta}{\delta u^\alpha}[P,Q]=D_P\left(\frac{\delta Q}{\delta u^\alpha}\right)
+(-1)^{pq}D_Q\left(\frac{\delta P}{\delta u^\alpha}\right);\nn\\
&   (-1)^{p-1}D_{[P,Q]}=D_P\circ D_Q-(-1)^{(p-1)(q-1)}D_Q\circ D_P.\nn
\end{align}
\end{lem}
\begin{prf}
The first identity can be obtained from the definition of $D_P$.
The second one is a corollary of the identity (iii) of Lemma 2.1.3 and the identity (i) of Lemma 2.1.5 given in  \cite{Jacobi}.
The third identity is a corollary of the second one. The lemma is proved.
\end{prf}

Theorem \ref{zh-01-22-a} gives the existence part of Theorem \ref{main-thm}, and Theorem \ref{thm-unq-1} (combining with
Theorem \ref{thm-unq-2}) gives the uniqueness part.

There are two important examples of such deformations when the flat exact semisimple bihamiltonian structures is provided
by a semisimple cohomological field theory. In \cite{DZ-NF} the first- and the third-named authors construct, for any semisimple Frobenius manifold,
the so-called topological deformation of the associated principal hierarchy and its tau structure. As we mentioned in Example \ref{zh-12-31b}, in \cite{Bu} Buryak
constructed a Hamiltonian integrable hierarchy associated to any cohomological field theory, and in \cite{BD} he and his collaborators
showed that this integrable hierarchy also possesses a tau structure. Buryak conjectured in \cite{Bu} that the above two integrable hierarchies
are equivalent via a Miura type transformation. He and his collaborators further refined this conjecture in \cite{BD} as an equivalence between tau-symmetric
Hamiltonian deformations via a normal Miura type transformation. The notion of  normal Miura type transformation was introduced in \cite{DLYZ}, our Definition \ref{equivalent} (see also Theorem \ref{thm-unq-2}) is a kind of its generalization. We hope our results
could be useful to solve the Buryak's \emph{et al} conjecture.

\section{Conclusion}\label{sec-7}

We consider in this paper the integrable hierarchies associated to a class of flat exact semisimple bihamiltonian structures of hydrodynamic type.
This property of flat exactness enables us to associate to any semisimple bihamiltonian structure of hydrodynamic type a Frobenius manifold
structure (without the Euler vector field), and a bihamiltonian integrable hierarchy which is called the principal hierarchy. We show that this
principal hierarchy possesses a tau structure and also the Galilean symmetry. For any deformation of the flat exact semisimple bihamiltonian
structures of hydrodynamic type which has constant central invariants, we construct the deformation of the principal hierarchy and show the
existence of tau structure and Galilean symmetry for this deformed integrable hierarchy. We also describe the ambiguity of the choice of tau
structure for the deformed integrable hierarchy. Our next step is to study properties of the Virasoro symmetries that are inherited from the
Galilean symmetry of the deformed integrable hierarchy in order to fix an appropriate representative of the tau structures which, in the case
associated to a cohomological field theory, corresponds to the partition function. We will do it in a subsequent publication.

\paragraph{Acknowledgements} This work  is partially supported by  NSFC No.\,11371214 and No.\,11471182.
B.D. kindly acknowledge the hospitality and generous support during his visit
to the Department of Mathematics of Tsinghua University where part of this work was completed.

\appendix

\section{On semi-Hamiltonian hierarchies}

In this appendix, we prove a classification theorem for conserved quantities and symmetries of a semi-Hamiltonian system satisfying certain nondegenerateness
conditions.

\begin{dfn}[\cite{tsarev}]\label{zh-12-31a}
i) A system of evolutionary partial differential equations of the form
\begin{equation}
\frac{\p u^i}{\p t}=A^i(u)u^{i,1}, \quad i=1, \dots, n \label{semi-0}
\end{equation}
is called semi-Hamiltonian, if $A^i\ne A^j$ for $i \ne j$ and
\[\p_k\left(\frac{\p_j A^i}{A^j-A^i}\right)=\p_j\left(\frac{\p_k A^i}{A^k-A^i}\right), \quad \mbox{for distinct }i, j, k, \]
where $\p_k=\frac{\p}{\p u^k}$.

ii) A semi-Hamiltonian system of evolutionary partial differential equations is called nondegenerate if $\p_i A^i \ne 0$ for all $i=1, \dots, n$.
We denote $A^i_i=\p_i A^i$ from now on.
\end{dfn}

According to Tsarev's results \cite{tsarev}, a semi-Hamiltonian system has infinitely many conserved quantities of the form
\[H=\int h \in\hF^0_0,\]
and infinitely many symmetries of the form
\[Y=\int\left(\sum_{i=1}^n B^i(u)u^{i,1}\theta_i\right)\in\hF^1_1.\]
An important question is: Are there conserved quantities and symmetries with higher degrees, which belong to 
$\hF^0_{\ge1}$ and $\hF^1_{\ge2}$? The following Theorem and Corollary give the answer.
\begin{thm}\label{app-thm}
Suppose $X$ is a nondegenerate semi-Hamiltonian system.
\begin{itemize}
\item[i)] If $H\in\hF^0_{\ge1}$ satisfies $[X, H]=0$, then $H=0$.
\item[ii)] If $Y\in\hF^1_{\ge2}$ satisfies $[X, Y]=0$, then $Y=0$.
\end{itemize}
\end{thm}
\begin{prf}
i) Suppose $H=\int h$, where $h\in \A^{(N)}$. Recall that $\A^{(N)}$ is the space of differential polynomials that do not depend on $u^{i,s}$ with $s>N$.
We denote $\delta_i=\frac{\delta}{\delta u^i}$ and define $Z_i=\delta_i[X,H]$. Since $[X, H]=0$, we know that $Z_i=0$. On the other hand, by using
Lemma \ref{last-lem} we have
\begin{align*}
Z_i &= D_X\left(\delta_i H\right)+D_H\left(\delta_i X\right)\\
&= \sum_{t\ge 0} \p^t\left(A^\alpha u^{\alpha,1}\right)\frac{\p \left(\delta_i H\right)}{\p u^{\alpha,t}}
+\p_i A^\alpha u^{\alpha, 1}\delta_\alpha H-\p\left(A^i \delta_i H\right).
\end{align*}
Then one can obtain that
\[0=(-1)^N \left(Z_i\right)_{(k, 2N+1)}=\left(A^k-A^i\right)h_{(k, N)(i, N)},\]
where $\left(f \right)_{(i, t)}=\frac{\p f }{\p u^{i,t}}$ for $f\in\A$.

The above equation implies that $h_{(k, N)(i, N)}=0$ for $i\ne k$, so we can assume that
\[h=\sum_{i} h_i(u, \dots, u^{(N-1)}; u^{i, N}).\]
If we can prove that $h_i$ depends on $u^{i, N}$ linearly, that is
\[h=\sum_{i} g_i(u, \dots, u^{(N-1)}) u^{i, N}+R(u, \dots, u^{(N-1)}),\]
then for $k\ne i$ we have
\[0=(-1)^N \left(Z_i\right)_{(k, 2N)}=\left(A^k-A^i\right)\left(\left(g_i\right)_{(k, N-1)}-\left(g_k\right)_{(i, N-1)}\right).\]
So there exists $g\in\A^{(N-1)}$ such that $g_{(i, N-1)}=g_i$. In particular, $h-\p g\in \A^{(N-1)}$. Then the first part of the
theorem can be proved by induction on $N$.

Now let us proceed to prove the linear dependence of $h_i$ on $u^{i, N}$. Define $Y_i=\left(h_i\right)_{(i, N)(i, N)}$, then we have the following identity:
\begin{align}
0=&(-1)^N \left(Z_i\right)_{(i, 2N)}\nn\\
=&\sum_{t\ge 0} \p^t\left(X^\alpha\right)\left(Y_i\right)_{(\alpha,t)}-A^i \p Y_i
+\left(2 A^i_i u^{i,1}+(2N-1) \p A^i\right)Y^i, \label{eq-Yi}
\end{align}
where $X^k=A^k u^{k,1}$. We need to show that $Y_i=0$. 

Consider $Y_i$ as a polynomial of $u^{i,1}, \dots, u^{i, N}$, and expand it as
\[Y_i=\sum_{\beta}c_{\beta} \left(u^{i,1}\right)^{\beta_1}\cdots\left(u^{i,N}\right)^{\beta_N},\]
where $\beta=(\beta_1, \dots, \beta_N)$, and $c_\beta$ do not depend on $u^{i,1}, \dots, u^{i, N}$. We define a lexicographical order on
the set of monomials recursively:
\begin{align*}
& \left(u^{i,1}\right)^{\beta_1}\cdots\left(u^{i,N}\right)^{\beta_N} \preceq \left(u^{i,1}\right)^{\gamma_1}\cdots\left(u^{i,N}\right)^{\gamma_N}\\
\Leftrightarrow \quad & \beta_N<\gamma_N \mbox{ or } \beta_N=\gamma_N \mbox{ and }  \\
& \left(u^{i,1}\right)^{\beta_1}\cdots\left(u^{i,N-1}\right)^{\beta_{N-1}}
\preceq \left(u^{i,1}\right)^{\gamma_1}\cdots\left(u^{i,N-1}\right)^{\gamma_{N-1}}.
\end{align*}
Then $Y_i$ can be written as the sum of its leading term and the remainder $Y_i=c\,w+R$, where
$w=\left(u^{i,1}\right)^{\beta_1}\cdots\left(u^{i,N}\right)^{\beta_N}$.

The equation \eqref{eq-Yi} now reads
\begin{align*}
0=& \left(\sum_{t \ge 0}\p^t X^\alpha (c)_{(\alpha,t)}-A^i \p c\right)\,w+\left(\sum_{t \ge 0}\p^t X^\alpha (w)_{(\alpha,t)}-A^i \p w\right)\,c\\
&+\left(2 A^i_i u^{i,1}+(2N-1)\p A^i\right)\,c\,w+\tilde{R},
\end{align*}
where $\tilde{R}$ is the sum of terms coming from $R$.

Consider the quotient of the above equation modulo $w$:
\begin{align}
0=&\sum_{t \ge 0}\p^t X^\alpha (c)_{(\alpha,t)}-A^i \p c+\left(\left(2+\beta_2+\cdots+\beta_N\right)A^i_i u^{i,1}\right.\nn\\
&
\left.+\left(2N-1+\beta_1+2\beta_2+\cdots+N \beta_N\right)\p A^i\right)c. \label{eq-c}
\end{align}
Suppose $c \in \A^{(m)}\ (m\ge 1)$, the derivative with respect to $u^{k,m+1}$ of the above equation gives
\[0=\left(A^k-A^i\right)(c)_{(k,m)} \quad \Rightarrow \quad (c)_{(k,m)}=0 \mbox{ for } k\ne i.\]
On the other hand $(c)_{(i,m)}=0$ by definition, so we have $c\in \A^{(m-1)}$. An induction on $m$ implies that $c\in \A_0$.

Finally, take the leading term of \eqref{eq-c} with respect to the lexicographical order, we obtain
\[\left(\beta_1+3\beta_2+\cdots+(N+1)\beta_N+2N+1\right)A^i_i u^{i,1}\,c=0,\]
so $c=0$, and consequently we deduce that $Y_i=0$. The first part of the theorem is proved.

ii) Suppose $Y=\int \left(B^\alpha\theta_\alpha\right)$, where $B^i\in \A^{(N)}$, and define $Z^i=\frac{\delta}{\delta \theta_i}[X, Y]$,
then we have
\[0=Z^i=\sum_{t\ge0}\p^t X^\alpha B^i_{(\alpha,t)}-B^\alpha X^i_{(\alpha,0)}-A^i\p B^i.\]
The equation $Z^i_{(k,N+1)}=0$ implies that $B^i_{(k,N)}=0$ for $k\ne i$. Denote $Y^i=B^i_{(i,N)}$, we have
\[0=Z^i_{i,N}=\sum_{t\ge0}\p^t X^\alpha Y^i_{(\alpha,t)}-A^i\p Y^i+N \p A^i\,Y^i-B^\alpha A^i_{(\alpha,0)}\delta_{N,1}.\]
When $N\ge2$, by using a similar argument that we used in the above proof of the first part of the theorem, we can show that the above equations have the only solution $Y^i=0$. When $N=1$,
suppose $B^i=c\, \left(u^{i,1}\right)^\beta+R$, then the leading term of the above equation reads
\[\left(\beta^2-1\right)A^i_i\,\left(u^{i,1}\right)^\beta c=0,\]
so we have $\beta=1$ or $c=0$. The second part of the theorem is proved.
\end{prf}

\begin{cor}\label{app-cor}
Given a vector field
\[X=X_1+X_2+\cdots\in\hF^1_{\ge1},\]
where $X_d\in \hF^1_d$, and $X_1$ is nondegenerate semi-Hamiltonian. Then the following statements hold true:
\begin{itemize}
\item[i)] If $H\in \F$ is a conserved quantity of $X$, then
\[H=H_0+H_1+\cdots\]
where $H_d\in\F_d$, $H_0$ is a conserved quantity of $X_1$, and $H_d\ (d\ge1)$ is uniquely determined by $H_0$.
\item[ii)] If $Y\in \hF^1_{\ge1}$ is a symmetry of $X$, then
\[Y=Y_1+Y_2+\cdots\]
where $Y_d\in\hF^1_d$, $Y_1$ is a symmetry of $X_1$, and $Y_d\ (d\ge2)$ is uniquely determined by $Y_1$.
\end{itemize}
\end{cor}

The proof is trivial, so we omit it. Note that the degree of a symmetry starts from $1$. There may exists other symmetries with degree zero
(e.g., the unit vector filed $Z$ for $X_{\alpha, 0}$).


\begin{thebibliography}{99}

\bibitem{BF} Behrend, K., Fantechi, B.: The intrinsic normal cone.
Invent. Math. \textbf{128}(1), 45--88 (1997)

\bibitem{Bu} Buryak A.: Double ramification cycles and integrable hierarchies, Comm. Math. Phys. \textbf{336}, 1085--1107  (2015)

\bibitem{BD} Buryak A., Dubrovin B., Gu\'er\'e J., Rossi P.: Tau-structure for the Double Ramification Hierarchies, eprint arXiv:\,1602.05423.

\bibitem{BD2} Buryak A., Dubrovin B., Gu\'er\'e J., Rossi P.: Integrable systems of double ramification type, eprint arXiv:\,1609.04059.

\bibitem{BPS-1} Buryak, A., Posthuma, H., Shadrin, S.:
On deformations of quasi-Miura transformations and the Dubrovin-Zhang bracket.
J. Geom. Phys. \textbf{62}(7), 1639--1651 (2012)

\bibitem{BPS-2} Buryak, A., Posthuma, H., Shadrin, S.:
A polynomial bracket for the Dubrovin-Zhang hierarchies.
J. Diff. Geom. \textbf{92}(1), 153--185 (2012)

\bibitem{CPS-1} Carlet, G., Posthuma, H., Shadrin, S.:
Bihamiltonian cohomology of KdV brackets.  Comm. Math. Phys. \textbf{341}(3), 805--819 (2016)

\bibitem{CPS-2} Carlet, G., Posthuma, H., Shadrin, S.:
Deformations of semisimple Poisson pencils of hydrodynamic type are unobstructed. eprint arXiv:\,1501.04295.

\bibitem{CPS-3} Carlet, G., Posthuma, H., Shadrin, S.:
The bi-Hamiltonian cohomology of a scalar Poisson pencil.
Bull. Lond. Math. Soc. \textbf{48}(4), 617--627 (2016)

\bibitem{DKJM} Date, E., Kashiwara, M., Jimbo, M., Miwa, T.: Transformation groups for soliton equations.
Nonlinear integrable systemsâclassical theory and quantum theory (Kyoto, 1981), 39--119, World Sci. Publishing, Singapore, 1983.

\bibitem{Di} Dickey, L. A.: Soliton equations and Hamiltonian systems.
Advanced Series in Mathematical Physics, 12. World Scientific Publishing Co., Inc., River Edge, NJ, 1991.

\bibitem{Du-1}  Dubrovin, B.: Geometry of 2D topological field theories.
Integrable systems and quantum groups (Montecatini Terme, 1993), 120--348, Lecture Notes in Math. \textbf{1620}, Springer, Berlin, 1996.

\bibitem{Du-2} Dubrovin, B.: Flat pencils of metrics and Frobenius manifolds.
Integrable systems and algebraic geometry (Kobe/Kyoto, 1997), 47--72, World Sci. Publ., River Edge, NJ, 1998.

\bibitem{Du-3} Dubrovin, B.: Painlev\'e transcendents in two-dimensional topological field theory. The Painlev\'e property, 287--412,
CRM Ser. Math. Phys., Springer, New York, 1999.

\bibitem{DZ-NF} Dubrovin, B., Zhang, Y.: Normal forms of hierarchies of integrable PDEs, Frobenius manifolds and Gromov--Witten invariants.
eprint  arXiv:\,math/0108160.

\bibitem{DZ-Toda}  Dubrovin, B., Zhang, Y.: Virasoro symmetries of the extended Toda hierarchy.
Comm. Math. Phys. \textbf{250}(1), 161--193 (2004)

\bibitem{DLZ-1} Dubrovin, B., Liu, S.-Q., Zhang, Y.: On Hamiltonian perturbations of hyperbolic systems of conservation
laws. I. Quasi-triviality of bi-Hamiltonian perturbations. Comm. Pure Appl.Math. \textbf{59}(4), 559--615 (2006)

\bibitem{DN83} Dubrovin, B., Novikov, S. P.:
Hamiltonian formalism of one-dimensional systems of the hydrodynamic type and the Bogolyubov-Whitham averaging method. (Russian)
Dokl. Akad. Nauk SSSR \textbf{270}(4), 781--785 (1983) 


\bibitem{DLYZ} Dubrovin, B., Liu, S.-Q., Yang D., Zhang, Y. : Hodge integrals and tau-symmetric integrable hierarchies of Hamiltonian evolutionary PDEs. Adv. Math. \textbf{293}, 382--435 (2016)

\bibitem{EF} Enriquez, B., Frenkel, E.: Equivalence of two approaches to integrable hierarchies of KdV
type. Comm. Math. Phys. \textbf{185}, 211--230  (1997)

\bibitem{FL} Falqui, G., Lorenzoni, P.: Exact Poisson pencils, $\tau$-structures and topological hierarchies.
Phys. D \textbf{241}(23-24), 2178--2187  (2012)

\bibitem{Fera} Ferapontov, E. V.: Compatible Poisson brackets of hydrodynamic type.
Kowalevski Workshop on Mathematical Methods of Regular Dynamics (Leeds, 2000).
J. Phys. A \textbf{34}(11), 2377--2388  (2001)

\bibitem{Ge}  Getzler, E.: The Toda conjecture.
Symplectic geometry and mirror symmetry (Seoul, 2000), 51--79, World Sci. Publ., River Edge, NJ, 2001.

\bibitem{JMU-1} Jimbo, M., Miwa, T., Ueno, K.: Monodromy preserving deformation of linear ordinary differential equations with rational coefficients. I.
Phys. D \textbf{2}(2), 306--352  (1981)

\bibitem{JM-2} Jimbo, M., Miwa, T.: Monodromy preserving deformation of linear ordinary differential equations with rational coefficients. II.
Phys. D \textbf{2}(3), 407--448  (1981)

\bibitem{JM-3}  Jimbo, M., Miwa, T.: Monodromy preserving deformation of linear ordinary differential equations with rational coefficients. III.
Phys. D \textbf{4}(1), 26--46 (1981/82)

\bibitem{KW} Kac, V., Wakimoto, M.: Exceptional hierarchies of soliton equations.
 Proc. Sympos. Pure Math. 49, part 1, 191--237 (1989) 

\bibitem{Ko} Kontsevich, M.: Intersection theory on the moduli space of curves and the matrix Airy function.
Comm. Math. Phys \textbf{147}(1), 1--23 (1992)

\bibitem{KM} Kontsevich, M., Manin, Yu.: Gromov--Witten classes, quantum cohomology, and enumerative geometry.
Comm. Math. Phys. \textbf{164}(3), 525--562 (1994)

\bibitem{LT} Li, J., Tian, G.: Virtual moduli cycles and Gromov--Witten invariants of algebraic varieties.
J. Amer. Math. Soc. \textbf{11}(1), 119--174 (1998)

\bibitem{LZ-1} Liu, S.-Q., Zhang, Y.: Deformations of semisimple bihamiltonian structures of hydrodynamic type.
J. Geom. Phys. \textbf{54}(4), 427â453 (2005)

\bibitem{Jacobi} Liu, S.-Q., Zhang, Y.: Jacobi structures of evolutionary partial differential equations.
Adv. Math. \textbf{227}(1), 73--130 (2011)

\bibitem{BCIH-I} Liu, S.-Q., Zhang Y.: Bihamiltonian Cohomologies and Integrable Hierarchies I: A Special Case.
Comm. Math. Phys. \textbf{324}, 897--935 (2013)

\bibitem{Mi-Toda}  Milanov, T.: The equivariant Gromov--Witten theory of $\mathbb{CP}^1$ and integrable hierarchies.
Int. Math. Res. Not. IMRN \textbf{2008}, Art. ID rnn 073, 21 pp.

\bibitem{Mira} Miramontes, J. L.: Tau-functions generating the conservation laws for generalized integrable hierarchies of KdV and affine Toda type. Nuclear Phys. \textbf{B 547}, 623â663 (1999)

\bibitem{OP} Okounkov, A., Pandharipande, R.: The equivariant Gromov--Witten theory of $\mathbb{P}^1$.
Ann. of Math. (2) \textbf{163}(2), 561--605 (2006)

\bibitem{RT}  Ruan, Y., Tian, G.: A mathematical theory of quantum cohomology.
J. Differential Geom. \textbf{42}(2), 259--367 (1995)

\bibitem{Sato} Sato, M.: Soliton equations as dynamical systems on infinite-dimensional Grassmann manifold. RIMS Kokyuroku \textbf{439} (1981) 30-46.

\bibitem{SS} Sato, M., Sato, Y.: Soliton equations as dynamical systems on infinite-dimensional Grassmann manifold. in: Nonlinear partial differential equations in applied science (Tokyo, 1982), 259â271, North-Holland Math. Stud., 81, North-Holland, Amsterdam, 1983. 

\bibitem{SW} Segal, G., Wilson, G.: Loop groups and equations of KdV type. 
Publ. Math. IHES 61, 5--65 (1985).

\bibitem{tsarev} Tsarev, S.: Geometry of Hamiltonian systems of hydrodynamic type. Generalized hodograph method. Izv. Akad. Nauk SSSR, Ser. Mat. (1990).

\bibitem{Wi} Witten, E.: Two-dimensional gravity and intersection theory on moduli space.
Surveys in differential geometry (Cambridge, MA, 1990) 1, Bethlehem, PA: Lehigh Univ., pp. 243--310.

\bibitem{Wu} Wu, C.-Z.: Tau Functions and Virasoro Symmetries for Drinfeld-Sokolov Hierarchies. eprint arXiv:1203.5750.

\bibitem{XZ} Xue, T., Zhang, Y.: Bihamiltonian systems of hydrodynamic type and reciprocal transformations. Lett. Math. Phys. \textbf{75}(1), 79--92 (2006)

\bibitem{yz} Zhang, Y.: Deformations of the bihamiltonian structures on the loop space of Frobenius manifolds. Recent advances in integrable systems (Kowloon, 2000). J. Nonlinear Math. Phys. \textbf{9} suppl. 1, 243--257 (2002)
\end{thebibliography}
\end{document}